\documentstyle[11pt]{article}
\begin{document}%
\input amssym.def
\input amssym
\font\tencmmib=cmmib10 \font\sevencmmib=cmmib7 \font\fivecmmib=cmmib5
\newfam\cmmibfam \textfont\cmmibfam=\tencmmib
\scriptfont\cmmibfam=\sevencmmib \scriptscriptfont\cmmibfam=\fivecmmib
\font\omegathic=eufm10 at 12 pt
\font\script=eusm10 at 12 pt

\def\bs{\backslash}

\def\ge{\epsilon}
\def\aa{\hbox{\omegathic a}}
\def\bb{\hbox{\omegathic b}} \def\cc{\hbox{\omegathic c}}
\def\pp{\hbox{\omegathic p}} \def\qq{\hbox{\omegathic q}}
\def\ll{\hbox{\omegathic l}} \def\Gg{\hbox{\omegathic g}}
\def\HH{\hbox{\omegathic H}} \def\scI{\hbox{\script I}}
\def\scC{\hbox{\script C}}
\newtheorem{theorem}{Theorem}[section] \newtheorem{claim}{Claim}[theorem]
\newtheorem{Lemma}{Lemma}[theorem]
\newtheorem{Corollary}{Corollary}[theorem]
\newtheorem{lemma}[theorem]{Lemma}
\newtheorem{proposition}[theorem]{Proposition}
\newtheorem{corollary}[theorem]{Corollary}
\newtheorem{observation}[theorem]{Observation}
\newtheorem{question}[theorem]{Question} \newcounter{def}
\newenvironment{definition}{\refstepcounter{theorem}
  \vskip.2cm\noindent{\bf Definition \thesection.\arabic{theorem}
    }}{\vskip.2cm\noindent} \newenvironment{fact}{\refstepcounter{theorem}
  \vskip.2cm\noindent{\bf Fact \thesection.\arabic{theorem}
    }}{\vskip.2cm\noindent }
\newenvironment{problem}{\refstepcounter{theorem} \vskip.2cm\noindent{\bf
    Problem \thesection.\arabic{theorem} }}{\vskip.2cm\noindent}
\newenvironment{notations}{\refstepcounter{theorem} \vskip.2cm\noindent{\bf
    Notations \thesection.\arabic{theorem} }}{\vskip.2cm\noindent}
\newenvironment{facts}{\refstepcounter{theorem} \vskip.2cm\noindent{\bf
    Facts \arabic{section}.\arabic{theorem} }}{\vskip.2cm\noindent}
\newenvironment{example}{\refstepcounter{theorem} \vskip.2cm\noindent{\bf
    Example \thesection.\arabic{theorem} }}{\vskip.2cm\noindent}
\newenvironment{condition}[1]{ \vskip.2cm\noindent{\bf Condition #1:
    }}{\vskip.2cm} \newenvironment{examples}{\refstepcounter{theorem}
  \vskip.2cm\noindent{\bf Examples \thesection.\arabic{theorem}:
    }}{\vskip.2cm\noindent}
\newenvironment{remark}{\refstepcounter{theorem} \vskip.2cm\noindent{\bf
    Remark \thesection.\arabic{theorem} }}{\vskip.2cm\noindent}
\newenvironment{remarks}{\refstepcounter{theorem} \vskip.2cm\noindent{\bf
    Remarks \thesection.\arabic{theorem} }}{\vskip.2cm\noindent}
\def\Ad{\hbox{Ad}} \def\ad{\hbox{ad}} \def\dim{\hbox{dim}}
\def\deg{\hbox{deg}} \def\Pic{\hbox{Pic}} \def\Jac{\hbox{Jac}}
\def\Aut{\hbox{Aut}} \def\Im{\hbox{Im}} \def\mod{\hbox{mod}}
\def\diag{\hbox{diag}} \def\det{\hbox{det}} \def\Ker{\hbox{Ker}}
\def\Hess{\hbox{Hess}} \def\rank{\hbox{rank}}
\def\hra{\hookrightarrow} \def\hla{\hookleftarrow}
\def\lra{\longrightarrow} \def\sura{\twoheadrightarrow}
\def\lla{\longleftarrow} \def\rar{\rightarrow}
\def\llra{\longleftrightarrow} \def\ra{\rightarrow} \def\Ra{\Rightarrow}
\def\La{\Leftarrow} \def\bs{\ifmmode {\setminus} \else$\bs$\fi}
\def\hat{\widehat} \def\tilde{\widetilde} \def\nin{\not\in}
\def\inn{\subset} \def\nni{\supset} \def\und{\underline }
\def\ove{\overline}
\def\sdprod{\rtimes}
\def\ende{\hfill $\Box$ \vskip0.25cm } \def\cD{\ifmmode {\cal D}
  \else$\cD$\fi}
\def\xg{\ifmmode {X_{\Gamma}} \else$\xg$\fi} \def\xgeq{\ifmmode {\xg =
    \Gamma\bs\cD} \else$\xgeq$\fi}

\def\xgs{\ifmmode {X_{\Gamma}^*} \else$\xgs$\fi}

\def\xgc{\ifmmode {\overline{X}_{\Gamma}} \else$\xgc$\fi}

\title{Hyperbolic planes} \author{Bruce Hunt \\ FB Mathematik,
  Universit\"at \\ Postfach 3049 \\ 67653 Kaiserslautern} \maketitle

\begin{center}
  {\Large \bf Introduction}
\end{center}

A {\it plane} is a two-dimensional {\it right} vector space $V$ over a
division algebra $D$, which we assume in this paper has an involution; a
hermitian form $h:V\times V \lra D$ is {\it hyperbolic}, if, with respect
to a properly choosen basis, it is given by the matrix ${0\ 1\choose 1\
  0}$; the pair $(V,h)$ is called a {\it hyperbolic plane}. If $D$ is
central simple over a number field $K$, a hyperbolic plane $(V,h)$ gives
rise to a reductive $K$-group $G_D=U(V,h)$ and a simple $K$-group
$SG_D=SU(V,h)$.  In this paper we study the case where the real Lie group
$G({\Bbb R})$ is of {\it hermitian type}, in other words that the symmetric
space of maximal compact subgroups of $G({\Bbb R})$ is a hermitian symmetric
space.  Let $d=\deg D$ (that is, $\dim_KD=d^2$); we consider the following
cases:
\begin{itemize}\item[1)] $d=1$: $k$ is a totally real number field of degree
  $f$ over ${\Bbb Q}$, $K|k$ is an imaginary quadratic extension, and $D:=K$,
  $V:=K^2$. The hyperbolic plane is denoted $(K^2,h)$.
\item[2)] $d=2$: $k$ as in 1), $D$ a totally indefinite quaternion division
  algebra, central simple over $k$, $V:=D^2$. The hyperbolic plane is
  denoted $(D^2,h)$.
\item[3)] $d\geq 3$: $k$ as above, $K|k$ an imaginary quadratic extension,
  $D$ a central simple division algebra over $K$ with a $K|k$-involution of
  the second kind, $V:=D^2$; this is again denoted $(D^2,h)$.
\end{itemize}
Of course, 1) is the same case as 3), if we view $K$ as central simple over
itself; the $K|k$ involution is just the Galois action of $K$ over $k$. Our
first observation (Proposition \ref{p5.1}) is that the case $d=1$ is none
other than the case familiar from the study of Hilbert modular varieties.
If $K=k(\sqrt{-\eta})$, with $\eta\in k^+$, then the isomorphism between
$SL_2(k)$ and $SU(K^2,h)$\footnote{In this paper we will reserve the
  notation $SL_2(k)$ for that group, and otherwise will use notations like
  $SL(n,K)$, etc.} is simply
\begin{eqnarray*} SL_2(k) & \stackrel{\sim}{\lra} & SU(K^2,h) \\
  \left(\begin{array}{cc} a & b \\ c & d \end{array}\right) & \mapsto &
  \left( \begin{array}{cc} d & 2c/\sqrt{-\eta} \\ b\sqrt{-\eta}/2 & a
\end{array}\right).\end{eqnarray*}
The following properties of Hilbert modular surfaces are well-known.
\begin{itemize}\item[i)] the domain $\cD\cong \HH\times \cdots \times \HH$
  is a product of $f$ copies of the upper half plane, and the rational
  boundary components are zero-dimensional, i.e., reduce to points.
\item[ii)] For the maximal order ${\cal O}_k\inn k$, the arithmetic subgroup
  $\Gamma:=SL_2({\cal O}_k)$ acts on $\cD$, the quotient
$X_{\Gamma}=\Gamma\bs \cD$ being
  non-compact, with a natural singular compactification $(\Gamma\bs \cD)^*$
  obtained by adding a point at each boundary component (the cusps).
  Furthermore, $X_{\Gamma}$ has exactly $h(k)$ cusps, where $h(k)$ denotes the
  class number of $k$.
\item[iii)] There are modular subvarieties, known in the case of Hilbert
  modular surfaces as Hirzebruch-Zagier cycles.
\end{itemize}
These are some of the properties it would be reasonable to consider for
hyperbolic planes over higher-dimensional $D$. The main results of this
paper contribute towards that end.  First, in analogy to i), we find that
for a hyperbolic plane $(D^2,h)$, the domain $\cD$ is a product of $f$
copies of an irreducible domain $\cD_1$, and $\cD_1$ is a tube
domain\footnote{this is certainly well-known to experts}. It is the disc
for $d=1$, Siegel space of degree two (type $\bf III_2$) for $d=2$, and of
type $\bf I_{d,d}$ for $d\geq3$, respectively. The rational boundary
components for the groups $G_D$ are points. It already follows from this
that the Baily-Borel compactification $\xgs$ of the arithmetic quotient
$\xg=\Gamma\bs \cD$, $\Gamma\inn G_D$ arithmetic, is obtained by adding
finitely many isolated
cusps to $\xg$, much as in the situation for Hilbert modular surfaces. We
remark that although we will not study the smooth compactifications of
$\overline{X}_{\Gamma}$, Satake in \cite{S} has considered these, and in
particular has given a formula for a ``cusp defect''.

To introduce arithmetic subgroups one fixes a maximal order $\Delta\inn D$,
and considers $\Gamma_{\Delta}=M_2(\Delta)\cap G_D$ and $S\Gamma_{\Delta}=
M_2(\Delta)\cap
SG_D$.  $S\Gamma_{\Delta}$ is the generalisation of the Hilbert modular group
($d=1$) to the higher $d$ case, and $\Gamma_{\Delta}$ is the generalisation of
the extended Hilbert modular group. The arithmetic quotients we study are
the $\Gamma_{\Delta}\bs \cD=X_{\Gamma_{\Delta}}$ and $S\Gamma_{\Delta}\bs
\cD=X_{S\Gamma_{\Delta}}$. The first surprise is that one can easily calculate
the class number, that is, the (finite) number of cusps on
$X_{\Gamma_{\Delta}}$.
We prove (Corollary \ref{c16.1})
\begin{theorem} The number of cusps of $\Gamma_{\Delta}$ is the class number of
  $K$ for $d\geq 3$ and the class number of $k$ for $d=1,2$.
\end{theorem}
This result is just a straightforward generalisation of the known $d=1$
result, as mentioned in ii). We then turn to a generalisation of iii)
above, and describe a rough analogue of the Hirzebruch-Zagier cycle $F_1$
to the case of higher $d$. This results from the enrichment of the
algebraic theory, which arises from {\it subfields} of the division algebra
$D$. Let us briefly explain this in the simplest case, $d\geq3$.  Then $D$
is a cyclic algebra $D=(L/K,\sigma,\gamma),\ \gamma\in K^*$, and the
subfield $L\inn
D$ is of degree $d$ over $K$; it is a splitting field for $D$. These
matters are recalled in the first paragraph for the convenience of the
reader, as they are an essential part of the remainder of the paper. The
subfield $L\inn D$ defines a subspace $L^2\inn D^2$, and $h_{|L^2}=:h_L$ is
again hyperbolic. In other words, in a hyperbolic plane over a degree $d$
cyclic division algebra, we have a subhyperbolic plane $(L^2,h)$; this is
of the type $d=1$, but with $f=d$.  This results in corresponding subgroups
$G_L\inn G_D$, subdomains $\cD_L\inn \cD$ and arithmetic subgroups
$\Gamma_{{\cal O}_L}\inn \Gamma_{\Delta}$; finally these yield modular
subvarieties
$M_L\inn X_{\Gamma_{\Delta}}$. By the result mentioned above,
this means: {\it on
  an arithmetic quotient $X_{\Gamma_{\Delta}}$ for $d\geq 3$, there are Hilbert
  modular subvarieties of dimension $d$}. For the general statement, see
Theorem \ref{t17.1} below. An additional analysis leads to the following
strengthening (Theorem \ref{t17a.1}):
\begin{theorem} Given a cusp
  $p\in X_{\Gamma_{\Delta}}^*\bs X_{\Gamma_{\Delta}}$ there
  is a modular subvariety $M_{L,p}\inn X_{\Gamma_{\Delta}}$, such that $p\in
  M_{L,p}^*$, that is, $p$ is also a cusp for $M_{L,p}$.
\end{theorem}
It is now quite straightforward to apply Shimura's theory to give the
moduli interpretation of the arithmetic quotients $X_{\Gamma_{\Delta}}$.
A second
application of this theory gives the moduli interpretation of the modular
subvarieties $M_L$. This is summarised as follows:
\begin{itemize}\item[1)] $d=1$: abelian varieties of dimension $2f$ with
  complex multiplication by $K$, of signature $(1,1)$, OR: abelian varieties
  of dimension $f$ with real multiplication by $k$.
\item[2)] $d=2$: abelian varieties of dimension $4f$, with multiplication
  by $D$.
\item[3)] $d=3$: abelian varieties of dimension $2d^2f$, with
  multiplication by $D$.
\end{itemize}
We note that the two possibilities for $d=1$ are related by the fact that,
given $A^f$, an abelian variety of dimension $f$ with real multiplication
by $k$, we get one of dimension $2f$ with complex multiplication as
follows: if the field $K=k(\sqrt{-\eta})$, let $K'={\Bbb Q}(\sqrt{-\eta})$.
Then
if $A^f \cong {\Bbb C}^f/\Lambda$, we have a lattice in ${\Bbb C}^{2f}$
given by
$\Lambda':=\Lambda\otimes_{{\Bbb Z}}{\cal O}_{K'}$.
The abelian variety $A^{2f}:= {\Bbb C}^{2f}/\Lambda'$ is
the desired abelian variety with complex multiplication. Conversely, if
$A^{2f}$ is given, it is easily seen that the lattice $\Lambda'$ has a
sublattice $\Lambda$ such that
$\Lambda'\cong \Lambda\otimes_{{\Bbb Z}}{\cal O}_{K'}$,
and ${\Bbb C}^f/\Lambda$
is an abelian variety of half the dimension with real multiplication.

The modular subvarieties, for $d\geq2$, parameterise:
\begin{itemize}\item[2)] $d=2$: the abelian variety $A^{4f}$ with
  multiplication by $D$ splits, $A^{4f}$ is isogenous to a product
  $A^{2f}\times A^{2f}$, where $A^{2f}$ has complex multiplication by an
  imaginary quadratic extension field $L\inn D$ of the center of $D$, $k$.
\item[3)] $d\geq 3$: the abelian variety $A^{2d^2f}$ with multiplication by
  $D$ splits, $A^{2d^2f}$ is isogenous to a product
  $\underbrace{A^{2df}\times \cdots \times A^{2df}}_{\hbox{$d$ factors}}$,
  where $A^{2df}$ has complex multiplication by the splitting field $L\inn
  D$.
\end{itemize}
In the fifth paragraph we discuss in more detail a central simple division
algebra $D$ of degree three over an imaginary quadratic extension
$K|{\Bbb Q}$, which defines unitary groups
which it appears may be related to Mumford's fake projective plane.

\section{Cyclic algebras with involution}
We will be considering matrix algebras $M_n(D)$, where $D$ is a simple
division algebra with an involution. There are two kinds of involution to
be considered. Fix, for the rest of the paper, a totally real number field
$k$ of finite degree $f$ over ${\Bbb Q}$, and assume $D$ is a division
algebra over $k$ with involution. The two cases are:
\begin{itemize}\item[(i)] $D$ is central simple over $k$, and $k$ is a
  maximal field which is symmetric with respect to the involution on $D$
  (involution of the first kind).
\item[(ii)] $D$ is central simple over an imaginary quadratic extension
  $K|k$, and $k$ is the maximal subfield of the center which is symmetric with
  respect to the involution on $D$; the involution restricts on $K$ to the
  $K|k$ involution (involution of the second kind).
\end{itemize}
If one sets $K:=k$ in case (i), then in both cases one speaks of a
$K|k$-involution. With these notations, let $d$ be the degree of $D$ over
$K$ (i.e., $\dim_KD=d^2$). Specifically, we will be considering the
following cases:
\begin{equation}\label{e1.0}\begin{minipage}{16cm} \begin{itemize}\item[1)]
      $d=1,\ D=K$ with an involution of the second kind, namely the $K|k$
      involution.
    \item[2)] $d=2,\ D$ a totally indefinite quaternion division algebra
      over $k$ with the canonical involution (involution of the first
      kind).
    \item[3)] $d\geq3,\ D$ a central simple division algebra of degree $d$
      over $K$ with a $K|k$-involution (involution of the second kind).
\end{itemize}\end{minipage}\end{equation}
Next recall that any division algebra over $K$ is a {\it cyclic algebra}
(\cite{A}, Thm.~9.21, 9.22). These algebras are constructed as follows. Let
$L$ be a cyclic extension of degree $d$ over $K$ and let $\sigma$ denote a
generator of the Galois group $Gal_{L/K}$. For any $\gamma\in K^*$, one
forms the algebra generated by $L$ and an element $e$,
\begin{equation}\label{e1.1}\begin{minipage}{15cm} $$A=(L/K,\sigma,\gamma):=
    L\oplus e L\oplus \cdots \oplus e^{d-1}L,$$
    $$e^d=\gamma,\hspace*{2cm} e\cdot z = z^{\sigma}\cdot e,\ \
    \forall_{z\in L}.$$
\end{minipage}\end{equation}
This algebra can be constructed as a subalgebra of $M_d(L)$, by setting:
\begin{equation}\label{e1.2} e=\left(\begin{array}{cccc} 0 & 1 & \cdots & 0
      \\ \vdots & & \ddots & \vdots \\ & & & 1 \\ \gamma & 0 & \cdots & 0
\end{array}\right),
\hspace*{2cm} z=\left(\begin{array}{cccc} z &\cdots &\cdots & 0 \\ &
    z^{\sigma} & & \\ & & \ddots & \\ 0 & \cdots & \cdots &
    z^{\sigma^{d-1}}\end{array}\right).
\end{equation}
One verifies easily that these matrices fulfill the relations (\ref{e1.1}),
and letting $B$ be the algebra in $M_d(L)$ generated by $e$ and $z\in L$,
we have $(L/K,\sigma,\gamma)\cong B$, and clearly $B\otimes_KL\cong
M_n(L)$, giving the explicit splitting. It is known that the cyclic algebra
$A=(L/K,\sigma,\gamma)$ is split if and only if $\gamma \in N_{L|K}(L^*)$
(\cite{A}, Thm.~5.14).  We will assume that $A$ is a division algebra; we
will denote it in the sequel by $D$.

Now assuming $d\geq2$, we consider the question of involutions on $D$. It
is well-known that for $d=2$, an involution of the second kind is the base
change of an involution of the first kind (\cite{A}, Thm.~10.21), so the
assumption above that for $d=2$ the involution is of the first kind is no
real restriction\footnote{perhaps somewhat more pragmatic is the remark
  that because of the exceptional isomorphism between domains of type ${\bf
    I_{2,2}}$ and ${\bf IV_4}$, the case $d=2$ and involution of the second
  kind can be translated into bilinear forms in eight variables, hence does
  not fit in well in the framework of hyperbolic planes.}.  Furthermore, if
$d\geq3$, there are no involutions of the first kind, since an involution
of the first kind gives an isomorphism $A\cong A^{op}$ ($A^{op}$ the
opposite algebra), implying that the class $[A]$ in the Brauer group $Br_K$
is of exponent one or two. Hence, for $d\geq3$ the restriction that the
involution is of the second kind is no restriction at all. Finally in case
$d=1$, an involution of the first kind on $D$ is trivial, so this case is
not hermitian, and we will not consider it.

Let us now consider the quaternion algebras $D$. Recall that $D$ is said to
be {\it totally definite} (respectively {\it totally indefinite}), if at
all real primes $\nu$, the local algebra $D_{\nu}$ is the skew field over
${\Bbb R}$ of the Hamiltonian quaternions ${\Bbb H}$ (respectively, if for
all real primes $\nu$, the local algebra $D_{\nu}$ is split). Recall also
that $D$ {\it ramifies} at a finite prime $\pp$ if $D_{\pp}$ is a division
algebra; the isomorphism class of $D$ is determined by the (finite) set of
primes $\pp$ at which it ramifies and its isomorphism class at those
primes (\cite{A}, Thm.~9.34). As a special case of cyclic algebras,
quaternion algebras can be displayed as algebras in $M_2(\ell)$, where
$\ell/k$ is a real quadratic extension, $\ell=k(\sqrt{a})$. Then there is
some $b\in k^*$ such that
$$e=\left(\begin{array}{cc} 0 & 1 \\ b & 0 \end{array}\right),
\hspace*{2cm} z=\left(\begin{array}{cc} z & 0 \\ 0 &
    z^{\sigma}\end{array}\right), \ \ z\in \ell,$$ where
$z^{\sigma}=(z_1+\sqrt{a}z_2)^{\sigma}=z_1-\sqrt{a}z_2$. In other words, we
can write for any $\alpha\in D$,
\begin{equation}\label{e2.1} \alpha=\left(\begin{array}{cc} a_0+a_1\sqrt{a} &
      a_2+a_3\sqrt{a} \\ b(a_2-a_3\sqrt{a}) &
      a_0-a_1\sqrt{a}\end{array}\right).
\end{equation}
Then the canonical involution is given by the involution
\begin{equation}\label{e2.2}
  \left(\begin{array}{cc} a & b \\ c & d \end{array}\right) \mapsto \left(
\begin{array}{cc} d & -b \\ -c & a\end{array}\right)
\end{equation}
on $M_2(\ell)$, and the norm and trace are just the determinant and trace
of the matrix
(\ref{e2.1}). We will use the notation $(a,b)$ to denote this algebra
$(\ell/k,\sigma,b)$, $\ell=k(\sqrt{a})$, if no confusion can arise from
this.  We remark also that for quaternion algebras we have
$$Tr_{D|k}(x)=x+\overline{x},\quad N_{D|k}(x)=x\overline{x},$$ where the
trace and norm are the reduced traces and norms. All traces and norms
occuring in this paper are the reduced ones unless stated otherwise.

In the case of involutions of the second kind, note first that the
$K|k$-conjugation extends to the splitting field $L$; its invariant
subfield $\ell$ is then a totally real extension of $k$, also cyclic with
Galois group generated by $\sigma$. We have the following diagram:
\begin{equation}\label{e3.1}
  \unitlength1cm
\begin{picture}(2,2)
  \put(.9,-.35){$k$} \put(-.2,1){$\ell$} \put(2,.7){$K$} \put(.9,1.85){$L$}
  \put(.1,1.1){\line(6,5){.8}} \put(.1,.9){\line(5,-6){.8}}
  \put(1.9,.65){\line(-6,-5){.8}} \put(1.9,.85){\line(-5,6){.8}}
\end{picture}
\end{equation}
and the conjugations on $L$ and $K$ give the action of the Galois group on
the extensions $L/\ell$ and
$K/k$; these are ordinary imaginary quadratic extensions.  There are
precise relations known under which $D$ admits a $K|k$-involution of the
second kind.
\begin{theorem}[\cite{A}, Thm.~10.18]\label{t4A.1} A cyclic algebra
  $D=(L/K,\sigma,\gamma)$ has an involution of the second kind $\iff$ there
  is an element $\omega\in \ell$ such that
  $$\gamma\overline{\gamma}=N_{K|k}(\gamma)=N_{\ell|k}(\omega)=\omega\cdot
  \omega^{\sigma} \cdots \omega^{\sigma^{d-1}}.$$
\end{theorem}
If this condition holds, then an involution is given explicitly by setting:
\begin{equation} \label{e4A.1} (e^k)^J=\omega\cdots
  \omega^{\sigma^{k-1}}(e^k)^{-1},
  \left(\sum e^iz_i\right)^J = \sum \overline{z}_i(e^i)^J,
\end{equation}
where $x\mapsto \overline{x}$ denotes the $L/\ell$-involution. In
particular for $x\in L$ we have
$$x^J=\overline{x}, \hbox{ and } x=x^J \iff x\in \ell.$$ Later it will be
convenient to have a description for when $x+x^J=0$. This results from the
following.
\begin{theorem}[\cite{A}, Thm.~10.10]\label{t4A.2} Given an involution $J$
  of the second kind on an algebra $A$, central simple of degree $d$ over
  $K$, there are elements $u_1,\ldots, u_d$, with $u_i=u_i^J$, such that
  $A$ is generated over $K$ by $u_1,\ldots, u_d$. Furthermore, there is an
  element $q\in A,\ q^J=-q,\ q^2\in k$, such that, as a $k$-vector space,
  $$A=A^++qA^+,$$ where $A^+=\{x\in A | x=x^J\}$.
\end{theorem}
If $x\in A$ is arbitrary, then ${1\over 2}(x+x^J)\in A^+$, while ${1\over
  2}(x-x^J)\in A^-$. For example, we have $e^i+(e^i)^J=:E^i\in A^+$,
and then we have an isomorphism
\begin{equation}\label{e4A.2} A^+\cong \ell\oplus E\ell \oplus \cdots
  \oplus E^{d-1}\ell.
\end{equation}
If, as above, $K=k(\sqrt{-\eta})$, $L=\ell(\sqrt{-\eta})$, then we may take
$\sqrt{-\eta}=q$ in the theorem above, and for elements in
\begin{equation}\label{e4A.3} qA^+\cong \sqrt{-\eta}\ell \oplus
  E\sqrt{-\eta}\ell\oplus \cdots \oplus E^{d-1}\sqrt{-\eta}\ell,
\end{equation}
we have $y=-y^J$. In particular, the dimension of
$qA^+$ is $d^2$ as a $k$-vector space, and the dimension of $A$ is $2d^2$.

\section{Algebraic groups}
Let $D$ be a central simple $K$-division algebra with a $K|k$-involution
$x\mapsto \overline{x}$ as discussed above (denoted $x\mapsto x^J$ above,
but unless it is a cause of confusion we fix the notation $x\mapsto
\overline{x}$ for the remainder of the paper), with $K|k$ imaginary
quadratic for $d=1,\geq3$ and $K=k$ for $d=2$. Consider the simple algebra
$M_2(D)$; it is endowed with an involution, $M\mapsto
{^t\overline{M}}=:M^*$. The group of units of $M_2(D)$ is denoted
$GL(2,D)$, its derived group is denoted by $SL(2,D)$. Consider a non-degenerate
hermitian form on $D^2$ given by
\begin{equation}\label{e4.1} h({\bf x},{\bf y})=
  x_1\overline{y}_2+x_2\overline{y}_1,
\end{equation}
where ${\bf x}=(x_1,x_2), {\bf y}=(y_1,y_2) \in D^2$. This means the
hermitian form is given by the matrix $H:={0\ 1 \choose 1\ 0}$. The unitary
and special unitary groups for $h$ are
\begin{equation}\label{e4.2} U(D^2,h)=\left\{g\in GL(2,D) \Big|
    gHg^*=H\right\}, \hspace*{1cm} SU(D^2,h)=U(D^2,h)\cap SL(2,D).
\end{equation}
The equations defining $U(D^2,h)$ are then
\begin{equation}\label{e4.3} U(D^2,h)=\left\{g=\left(\begin{array}{cc}a & b
        \\ c & d \end{array}\right)\Big| a\overline{d}+b\overline{c}=1,
    a\overline{b}+b\overline{a}=c\overline{d}+d\overline{c}=0\right\}.
\end{equation}
The additional equation defining $SU(D^2,h)$ can be written in terms of
determinants, using Dieudonn\'e's theory of determinants over skew fields,
(see \cite{Ar}, p.~157)
\begin{equation}\label{e4.4} \det(g)=N_{D|k}(ad-aca^{-1}b)=1.
\end{equation}
The center of $U(D^2,h)$, which we will denote by ${\cal C}$, is given by
$${\cal C}=\left\{\left(\begin{array}{cc} a & 0 \\ 0 & a\end{array}\right)
  \Big| a\in K, a\overline{a}=1\right\}\cong U(1)\cap K,$$ where we view
$K\inn {\Bbb C}$ as a subfield of the complex numbers. In particular:
\begin{lemma}\label{l4.1} For $d=2$, $K=k$ a real field, we have
  ${\cal C}=\{\pm1\}$.
\end{lemma}
In other words, for the case $d=2$, there is no essential difference
between the unitary and special unitary groups. Otherwise, $U(D^2,h)$ is a
reductive $K$-group, and $SU(D^2,h)$ is the corresponding simple group.
\begin{proposition}\label{p4.1} Let $SG_D$ denote $SU(D^2,h)$, the simple
  algebraic $K$-group defined by the relations (\ref{e4.3}) and
  (\ref{e4.4}).  Then $SG_D$ is simple with the following index (cf.
  \cite{tits})
\begin{itemize}\item[1)] $d=1$, ${^2A_{1,1}^1}$.
\item[2)] $d=2$, $C_{2,1}^2$.
\item[3)] $d\geq 3$, ${^2A_{2d-1,1}^d}$.
\end{itemize}
\end{proposition}
{\bf Proof:} This is immediate from the description in \cite{tits} of these
indices. \ende The unitary group for $d=1$ is familiar in different terms.
If we set $A=\left(\begin{array}{cc}-1/2 & -1 \\ -1/2 &
    1\end{array}\right)$, then
$$A\left(\begin{array}{cc}0 & 1 \\ 1 &
    0\end{array}\right){^tA}=\left(\begin{array}{cc} 1 & 0 \\ 0 & -1
\end{array}\right).$$
In particular, these two hermitian forms are equivalent. For the form
$H_1:={1\ \ 0 \choose \ 0\ -1}$, the corresponding unitary group is
customarily denoted by $U(1,1)$. Hence for $K\inn {\Bbb C}$ imaginary
quadratic over a real field, we have an isomorphism
\begin{eqnarray}\label{e4.5} U(K^2,h) & \stackrel{\sim}{\lra} & U(1,1;K) \\
  g & \mapsto & A g A^{-1}; \nonumber
\end{eqnarray}
the group $U(1,1;K)$ is given by the familiar conditions
$$U(1,1;K)=\left\{ \left(\begin{array}{cc} \alpha c & \beta c \\
      \overline{\beta}\overline{c} &
      \overline{\alpha}\overline{c}\end{array}\right) \Big| \alpha,\beta\in
  K, \alpha\neq 0, c\overline{c}=1,
  \alpha\overline{\alpha}-\beta\overline{\beta}=1\right\},$$ and the
special unitary group is given by
$$SU(1,1;K)=\left\{\left(\begin{array}{cc} \alpha & \beta \\
      \overline{\beta} & \overline{\alpha}
\end{array}\right)\in GL(2,K) \Big|
\alpha\overline{\alpha}-\beta\overline{\beta}=1\right\}.$$
Recall that this latter group is isomorphic to a subgroup of $GL(2,{\Bbb
  R})$, in fact we have the well-known isomorphism
\begin{eqnarray}\label{e5.1} SU(1,1;{\Bbb C}) & \stackrel{\sim}{\lra} &
  SL(2,{\Bbb R}) \\ g & \mapsto & B g B^{-1}, \nonumber
\end{eqnarray}
where $B={-i\ \ i \choose \ \! 1\ \ 1}$ is the matrix of a fractional
linear transformation mapping the disk to the upper half plane. If, instead
of $B$, we use
\begin{equation}\label{e5.2} B_{\eta}=\left(\begin{array}{cc} -i & i \\
      \sqrt{\eta} & \sqrt{\eta} \end{array}\right),\hspace*{.75cm}
  K=k(\sqrt{-\eta}),
\end{equation}
then we have in fact
\begin{proposition}\label{p5.1} The matrix $B_{\eta}$ gives an isomorphism
  $$\begin{array}{rcl} SU(1,1;K) & \stackrel{\sim}{\lra} & SL_2(k) \\
    \mapsto & B_{\eta} g B_{\eta}^{-1}. \end{array}$$ In particular, we
  have an isomorphism
  $$\begin{array}{rcl} SU(K^2,h) & \stackrel{\sim}{\lra} & SL_2(k) \\ g &
    \mapsto & B_{\eta}A g A^{-1} B_{\eta}^{-1}.\end{array}$$
\end{proposition}
{\bf Proof:} We only have to check that $B_{\eta}gB_{\eta}^{-1}$ is real,
as it is clearly a matrix in $SL(2,K)$, and $SL(2,K)\cap SL(2,{\Bbb
  R})=SL_2(k)$.  If $g={\alpha\ \beta \choose \overline{\beta}\
  \overline{\alpha}}$ then
$$B_{\eta}gB_{\eta}^{-1} = {1 \over -2\sqrt{-\eta}}\left(\begin{array}{cc}
    \sqrt{-\eta}(-\alpha-\overline{\alpha}+\beta+\overline{\beta}) &
    -\alpha-\beta+\overline{\alpha}+\overline{\beta} \\
    \eta(\alpha-\beta+\overline{\beta}-\overline{\alpha}) &
    -\sqrt{-\eta}(\alpha+\beta+\overline{\alpha}+\overline{\beta})
    \end{array}\right)
=$$
$$\left(\begin{array}{cc} Re(\alpha)-Re(\beta) & Im(\alpha)+Im(\beta) \\
    \eta(-Im(\alpha)+Im(\beta)) & Re(\alpha)+Re(\beta)
  \end{array}\right),$$ which is clearly real. \ende Note that the inverse
$SL_2(k)\lra SU(K^2,h)$ is given by conjugating with
$C=(B_{\eta}A)^{-1}=\left(\begin{array}{cc} 0 & -2i \\ \sqrt{\eta} & 0
\end{array}\right)^{-1} = \left(\begin{array}{cc}0 & {1\over \sqrt{\eta}}
  \\ {i\over 2} & 0 \end{array}\right)$, so we may describe the latter
group as:
\begin{equation}\label{e5.3} SU(K^2,h)=\left\{\left(\begin{array}{cc} \delta &
        2\gamma/\sqrt{-\eta} \\ \beta\sqrt{-\eta}/2 & \alpha
      \end{array}\right) \Big| \alpha,\beta,\gamma,\delta\in k,
    \alpha\delta-\gamma\beta=1\right\}.
\end{equation}
We might remark at this point that it is a bit messy to describe the
unitary groups in the same fashion, as they are {\it not} isomorphic to
subgroups of $GL(2,{\Bbb R})$. We have the following extensions:
\begin{equation}\label{e5.4}\begin{array}{ccccccccc} 1 & \lra & SU(1,1;K) &
    \lra & U(1,1;K) & \lra & U(1)\cap K & \lra 1 \\ & & |\wr \\ 1 & \lra &
    SL_2(k) & \lra & GL(2,k) & \lra & k^* & \lra 1, \end{array}
\end{equation}
the first being an extension by a compact torus, the second by a split
torus. That is why the map of Proposition \ref{p5.1}, if extended to
$U(1,1;K)$, does not land in $GL(2,{\Bbb R})$. It is also appropriate to
remark here that, for $d=2$, the {\it reductive} group associated with the
problem is the group of symplectic similtudes, i.e., $g\in GL(2,D)$ such
that $gHg^*=H\lambda$ for some non-zero $\lambda\in D$. This gives an
extension of $SU(D^2,h)$ similar to the second sequence above.

Next we describe some subgroups of $G_D=U(D^2,h)$ and $SG_D$.  For this,
recall that we are dealing with right vector spaces, and matrix
multiplication is done from the right, i.e.,
$$(x_1,x_2)\left(\begin{array}{cc}a & b \\ c & d
\end{array}\right)=(x_1a+x_2c,
x_1b+x_2d).$$ First, the vectors $(1,0)$ and $(0,1)$ are isotropic; the
stabilisers of the lines they span in $D^2$ are maximal (and minimal)
parabolics. For example, the stabiliser of $(0,1)D$ is
\begin{equation}\label{e6.1} P=\left\{g \in G_D \Big| (0,1)g=(0,1)\lambda,
    \lambda\in D^*\right\} = \left\{g=\left(\begin{array}{cc}a & b\\ 0 &
        \overline{a}^{-1}\end{array}\right) \Big| a\in D^*,\
    b+\overline{b}=0 \right\}.
\end{equation}
Indeed, it is immediate that for $g={a\ b\choose c\ d}\in P$, we must have
$c=0$, while the first relation of (\ref{e4.3}) reduces to
$a\overline{d}=1$, or $d=\overline{a}^{-1}$. Since a parabolic is a
semidirect product of a Levi factor and the unipotent radical, to calculate
the latter we may assume $a=1$.  Then the second relation of (\ref{e4.3})
is $a\overline{b}+b\overline{a}=\overline{b}+b=0$. The corresponding
parabolic in $SG_D$ takes the form
\begin{equation}\label{ePar1} SP=\left\{g\in SG_D \Big|
    g=\left(\begin{array}{cc} a & b \\ 0 & \overline{a}^{-1}
      \end{array}\right), \ a\in D^*,\ b+\overline{b}=0 \right\}.
\end{equation}
Note that $g$ in (\ref{ePar1}) has $K$-determinant
$$\det_K(g)=N_{D|K}(a
\overline{a}^{-1})=N_{D|K}(a)N_{D|K}(\overline{a})^{-1},$$ so the
requirement on $g$ is that $a$ must satisfy
$N_{K|k}(N_{D|K}(a)/N_{D|K}(\overline{a}))=1$. This is satisfied for all
$a\in K^*$, where $a$ is viewed as an element of $D^*$; if $a\not\in K$,
then the condition becomes $N_{D|K}(a)=\overline{N_{D|K}(a)}$, i.e.,
$N_{D|K}(a)\in k$. It follows that the Levi component is a product
$$L=\left\{g=\left(\begin{array}{cc}a & 0 \\ 0 &
      \overline{a}^{-1}\end{array}\right) \Big| a\in D^*\right\}\cong
\left(U(1)\cap K\right)\cdot \left\{ g = \left(\begin{array}{cc}a & 0 \\ 0
      & \overline{a}^{-1}\end{array}\right) \Big| a\in D^*,\ N_{D|K}(a)\in
  k\right\}.$$ Now recall from (\ref{e4A.3}) that for $d\geq3$, the set of
elements fulfilling $b+\overline{b}=0$ in $D$ is $\sqrt{-\eta}D^+\cong
\sqrt{-\eta}\ell \oplus E\sqrt{-\eta}\ell \oplus \cdots \oplus
E^{d-1}\sqrt{-\eta}\ell$, where $E^i=e^i+\overline{e^i}$. In particular,
the unipotent radical of $SP$, generated by ${1\ b\choose 0\ 1}$ with
$b+\overline{b}=0$, has dimension $d^2$ over $k$.

We also have non-isotropic vectors $(1,1)$ and $(1,-1)$. Any element $g\in
G_D$ preserving one preserves also the other (as they are orthogonal).
Hence the stabiliser of $(1,1)$ is
\begin{equation}\label{e6.2}\label{eCom}
  {\cal K}=\left\{g=\left(\begin{array}{cc}
        a & c \\ c & a\end{array}\right) \Big|
    a\overline{a}+c\overline{c}=1,\ a\overline{c}+c\overline{a}=0\right\}.
\end{equation}
For $d=1$, this can be more precisely described as
$${\cal K}=\left\{ g=\left(\begin{array}{cc} a & {c \over \sqrt{-\eta}} \\
      -{c\sqrt{-\eta}\over \eta} & a \end{array}\right),\ a,c\in k,
  a^2-{c^2\over \eta}=1\right\}.$$ For $d=2$, the relations give, viewed as
equations over $k$, four relations, so the group ${\cal K}$ is
four-dimensional. Finally, for $d\geq3$, we get (over $k$) $d$ relations
(as elements $a\overline{a}$ are in $\ell$, so diagonal) from the first
condition and $d^2-d+1$ conditions from the second, leaving $d^2+d+2$
parameters for the group ${\cal K}$.

If $D$ has an involution of the second kind, we have the splitting field
$L$ ($=K$ for $d=1$), which is an imaginary quadratic extension of the
totally real field $\ell$ ($=k$ for $d=1$), see the diagram (\ref{e3.1}).
Suppose then $d=2$, $D$ quaternionic, say $D=(\ell/k,\sigma,b)$, where
$\ell=k(\sqrt{a}),\ e^2=b\in k^*$; in this case $\ell$ is a splitting
field. Recall the conjugation $\sigma$ is given by
$$(z_1+\sqrt{a}z_2)^{\sigma}=z_1-\sqrt{a}z_2,$$ so the element
$c=diag(\sqrt{a},-\sqrt{a})$ representing $\sqrt{a}\in \ell$ satisfies
$c^{\sigma}=-c$, while $e^{\sigma}=e$. Consequently, the relation
(\ref{e1.1}) for $c$ is $ec=c^{\sigma}e=-ce$, and
\begin{equation}\label{e7.1} (ec)^2=(ec)(-ce)=-ec^2e=-ab,
\end{equation}
so $k(ec)\cong k(\sqrt{-ab})$. If we assume, as we may, that $a>0, b>0$
(otherwise replace $-ab$ in what follows by the negative one of $a, b$),
then $L:=\sqrt{-ab}$ is an imaginary quadratic extension of $k$ which is a
subfield of $D$. So in all cases ($d=1,2,\geq3$) we have an imaginary
quadratic extension of $\ell$ ($d=1,\geq3$) or $k$ ($d=2$), $L$, which is a
subfield of $D$. We now claim that $U(L^2,h)$ is a natural subgroup of
$U(D^2,h)$. In all cases $U(L^2,h)$ and $U(D^2,h)$ are defined by the same
relations (\ref{e4.3}), so we have
\begin{equation}\label{e7.2} U(L^2,h)=U(D^2,h)\cap M_2(L),
\end{equation}
which verifies the claim. In the case $d=1$, where we took $L=K$, this is
in fact the whole group. But even in this case we can get subgroups in this
manner. Consider taking, instead of $L=K$, $L={\Bbb Q}(\sqrt{-\eta})$ for
$K=k(\sqrt{-\eta})$. More generally, for any ${\Bbb Q}\subseteq k'
\subseteq k$, we have the subfield $k'(\sqrt{-\eta})\inn K$, and we may
consider the corresponding subgroup:
$$U(k'(\sqrt{-\eta})^2,h) = U(K^2,h)\cap M_2(k'(\sqrt{-\eta})).$$ So,
setting $L=k'(\sqrt{-\eta})$ in (\ref{e7.2}), these subgroups are defined
by the same relations (\ref{e4.3}). Now recall the description (\ref{e5.3})
for groups of the type (\ref{e7.2}).
\begin{proposition}\label{p7.1} We have the following subgroups of
  $U(D^2,h)$.
\begin{itemize}\item[1)] $d=1$, for any $ {\Bbb Q}\subseteq k' \subseteq k$ we
  have the subgroup
  $$SL_2(k')\cong \left\{ \left(\begin{array}{cc} \alpha &
        2\beta/\sqrt{-\eta} \\ \gamma\sqrt{-\eta}/2 & \delta
      \end{array}\right) \Big| \alpha,\beta,\gamma,\delta \in k',
    \alpha\delta-\gamma\beta=1\right\}\inn U(K^2,h).$$
\item[2)] $d=2$, for the field $L=k(\sqrt{-ab})$ above, we have the
  subgroup
  $$SL_2(k)\cong \left\{ \left(\begin{array}{cc} \alpha & 2\beta/\sqrt{-ab}
        \\ \gamma\sqrt{-ab}/2 & \delta \end{array}\right) \Big|
    \alpha,\beta,\gamma,\delta \in k,
    \alpha\delta-\gamma\beta=1\right\}\inn U(D^2,h).$$
\item[3)] $d\geq3$, for the degree $d$ cyclic extension $L/K$ we have the
  following subgroup
  $$SL_2(\ell)\cong \left\{ \left(\begin{array}{cc} \alpha &
        2\beta/\sqrt{-\eta} \\ \gamma\sqrt{-\eta}/2 & \delta
      \end{array}\right) \Big| \alpha,\beta,\gamma,\delta \in \ell,
    \alpha\delta-\gamma\beta=1\right\}\inn U(D^2,h).$$
\end{itemize}
\end{proposition}
These are interesting subgroups, whose existence is a natural part of the
description of $D$ as a cyclic algebra. In all cases we may view these
subgroups as stabilisers. Let $L$ denote one of the fields
$k'(\sqrt{-\eta}), k(\sqrt{-ab}), \ell(\sqrt{-\eta})$ as in Proposition
\ref{p7.1}; let $L^2\inn D^2$ denote the $k$-vector subspace of $D^2$,
viewing $D$ as a $k$-vector space, and consider the stabiliser in $G_D$.
Noting that $L^2\inn D^2$ may be spanned over $L$ by two non-isotropic
vectors $(1,1)$ and $(1,-1)$, it is clear that $g\in U(L^2,h)$ preserves
the $L$ span of these two vectors, that is, $L^2$, hence $U(L^2,h)$ is
contained in the stabiliser. The converse is also true: if an endomorphism
of $D^2$ preserves the $L$ span, then its matrix representation is in
$M_2(L)$, so the stabiliser is contained in the intersection $U(D^2,h)\cap
M_2(L)=U(L^2,h)$. Let us record this as
\begin{observation} The subgroups of Proposition \ref{p7.1} are stabilisers
  (with determinant 1) of $k$-vector subspaces of the $k$-vector space
  $D^2$.
\end{observation}

\section{Domains}
Recall that given an almost simple algebraic group $G$ over $k$, by taking
the ${\Bbb R}$ points one gets a semisimple real Lie group $G({\Bbb R})$.
For the three kinds of groups above we determine these now.
\begin{proposition}\label{p9.1} The real groups $G_D({\Bbb R})$ are the
  following.
\begin{itemize}\item[1)] $d=1$, $G_D({\Bbb R})=(U(1,1))^f, SG_D({\Bbb R})\cong
  (SL(2,{\Bbb R}))^f$.
\item[2)] $d=2$, $G_D({\Bbb R})=SG_D({\Bbb R})=(Sp(4,{\Bbb R}))^f$.
\item[3)] $d\geq3$, $G_D({\Bbb R})=(U(d,d))^f,\ SG_D({\Bbb
    R})=(SU(d,d))^f$.
\end{itemize}\end{proposition}
{\bf Proof:} The $d=1$ case is obvious. For $d=2$, note that since $D$ is
totally indefinite, $D_{\nu}\cong M_2({\Bbb R})$ for all real primes, and
consequently $(G_D)_{\nu}$ is a simple group of type $Sp$ for each $\nu$,
not compact as $G_D$ is isotropic, of rank 2, hence $Sp(4,{\Bbb R})$.
Suppose now that $d\geq3$. By Proposition \ref{p4.1} we know the index of
$G_D$ is ${^2A_{2d-1,1}^d}$, and the index, as displayed in \cite{tits}, is
\begin{equation}\label{e9.1}
\setlength{\unitlength}{0.006500in}%
\begin{picture}(894,94)(53,693)
  \thicklines \put(140,780){\circle*{14}} \put( 60,700){\circle*{14}}
  \put(140,700){\circle*{14}} \put(400,700){\circle*{14}}
  \put(400,780){\circle*{14}} \put(460,740){\circle{14}}
  \put(600,700){\circle*{14}} \put(600,780){\circle*{14}}
  \put(680,700){\circle*{14}} \put(680,780){\circle*{14}}
  \put(860,780){\circle*{14}} \put(860,700){\circle*{14}}
  \put(940,780){\circle{14}} \put(940,700){\circle{14}} \put(
  65,700){\line( 1, 0){130}} \put( 60,780){\circle*{14}} \put(
  65,780){\line( 1, 0){130}} \put(940,770){\line( 0,-1){ 65}}
  \put(330,780){\line( 1, 0){ 70}} \put(330,700){\line( 1, 0){ 70}}
  \put(405,780){\line( 3,-2){ 45}}
  \put(400,700){\line( 5, 3){ 50}}
  \put(605,700){\line( 1, 0){ 70}} \put(605,780){\line( 1, 0){ 75}}
  \put(680,780){\line( 1, 0){ 50}} \put(680,700){\line( 1, 0){ 45}}
  \put(800,780){\line( 1, 0){ 60}} \put(800,700){\line( 1, 0){ 60}}
  \put(865,780){\line( 1, 0){ 70}} \put(860,700){\line( 1, 0){ 75}}
  \put(230,695){$\cdots\cdots$} \put(750,695){$\cdots$}
  \put(230,775){$\cdots\cdots$} \put(750,775){$\cdots$} \put(230,680){$d$
    even} \put(750,680){$d$ odd}
\end{picture}
\end{equation}

\noindent where the number of black vertices is $2d-2$,
$d-1$ on each branch. This
shows in particular that the {\it isotropic} root is the one farthest to
the right. If a root is isotropic over ${\Bbb Q}$, then all the more over
${\Bbb R}$.  Hence, considering the Satake diagram of the real groups of
type ${^2A}$, we see that the only possibility is $SU(d,d)$, as this is the
only real group of this type for which the right most vertex is ${\Bbb
  R}$-isotropic. Actually, this only proves that {\it at least one} factor
$(G_D)_{\nu}$ is $SU(d,d)$; but our hyperbolic form has signature $(d,d)$
at {\it any} real prime, so it holds for {\it all} (non-compact) factors.
That there are no compact factors follows from the fact that $G_D$ is
isotropic. In all cases there are $f$ factors, as this is the degree of
$k|{\Bbb Q}$, so the ${\Bbb Q}$-group, $Res_{k|{\Bbb Q}}G_D$, has $f$
factors over ${\Bbb R}$. This verifies the proposition in all cases. \ende
It is a consequence of this proposition that the symmetric spaces
associated to the real groups $G_D({\Bbb R})$ are {\it hermitian
  symmetric}; indeed, this is why we have choosen those $D$ to deal with.
More precisely, one has the following statement.
\begin{theorem}\label{t9.1} The hermitian symmetric domains defined by the
  $G_D({\Bbb R})$ as in Proposition \ref{p9.1}
\begin{itemize}\item are tube domains, products of irreducible components
  of the types $\bf I_{1,1}$, $\bf III_2$, $\bf I_{d,d}$ in the cases
  $d=1$, $d=2$ and $d\geq3$, respectively.
\item For any rational parabolic $P\inn G_D$, the corresponding boundary
  component $F$ such that $P({\Bbb R})=N(F)$ is a point.
\end{itemize}
\end{theorem}
{\bf Proof:} The first statement follows immediately from Proposition
\ref{p9.1}. The second is well known in the case $d=1$ from the study of
Hilbert modular varieties. For $d=2$ the statement clearly holds for
$k={\Bbb Q}$ (just look at the index, from which one sees the boundary
components are points, not one-dimensional), and the general statement
follows from this: all boundary components are of the type
$(c_1,\ldots,c_f)\in \cD_1^*\times \cdots \times \cD_f^*,$ where $ c_i\in
\cD_i^*$ is a point, which is a general fact about ${\Bbb Q}$-parabolics in
a {\it simple} ${\Bbb Q}$-group (this applies directly to $Res_{k|{\Bbb
    Q}}SG_D$, but $G_D$ and $SG_D$ define the same domain and boundary
components).  Finally, for $d\geq3$, we point out that in the diagram
(\ref{e9.1}), the isotropic vertex (farthest right) corresponds to a
maximal parabolic which stabilises a point, establishing the statement for
$k={\Bbb Q}$. Then the general statement follows as above from this. \ende

We now consider the structure of the real parabolics for $d\geq2$ in more
detail. For this we refer to Satake's book, \S III.4.
\begin{itemize}\item $d=2$: Here we have the case denote $b=n$ in Satake's
  book, and the parabolic in $Sp(4,{\Bbb R})$ is
  $$P({\Bbb R})=(G^{(1)}\cdot G^{(2)})\sdprod U\cdot V = 1\cdot GL(2,{\Bbb
    R})\sdprod Sym_2({\Bbb R}),$$ where $Sym_2({\Bbb R})$ denotes the set
  of symmetric real $2\times 2$ matrices. Refering to (\ref{e6.1}), we may
  identify the Levi component of $P$ with
  $$L\cong \left\{\left(\begin{array}{cc}a & 0 \\ 0 &
        \overline{a}^{-1}\end{array}\right) \Big| a\in D^*\right\} \cong
  D^*,$$ and, as $D$ is indefinite, $D^*({\Bbb R})\cong GL(2,{\Bbb R})$;
  this is the factor $GL(2,{\Bbb R})$. The unipotent radical is
  $$U\cong \left\{\left(\begin{array}{cc} 1 & b \\ 0 & 1\end{array}\right)
    \Big| \ Tr(b)=0\right\} \cong D^0 \hbox{ (=totally imaginary
    elements)},$$ which is three-dimensional, and an isomorphism $U({\Bbb
    R})\cong Sym_2({\Bbb R})$ is given by
  $$\left(\begin{array}{cc}a_1\sqrt{a} & a_2+a_3\sqrt{a} \\
      b(a_2-a_3\sqrt{a}) & -a_1\sqrt{a}\end{array}\right) \mapsto
  \left(\begin{array}{cc} a_1 & a_2 \\ a_2 & a_3\end{array}\right).$$ This
  describes the real parabolic completely in this case.
\item $d\geq3$: Here we have $b=d$ in Satake's notation, $d=p=q$. The real
  parabolic is $$P({\Bbb R})=(G^{(1)}\cdot G^{(2)})\sdprod U\cdot V =
  (U(1)\cdot GL^0(d,{\Bbb C}))\sdprod {\cal H}_d({\Bbb C}),$$ where
  $GL^0(d,{\Bbb C})=\{g\in GL(d,{\Bbb C})|\det(g)\in {\Bbb R}^*\}$ and
  ${\cal H}_d({\Bbb C})$ denotes the space of hermitian $d\times d$
  matrices. We must now consider the parabolic $SP$ as in (\ref{ePar1}).
  The Levi component is again
  $$L\cong \left\{\left(\begin{array}{cc} a & 0 \\ 0 & \overline{a}^{-1}
    \end{array}\right) \Big| a\in D^*,\
  N_{D|k}(a\overline{a}^{-1})=1\right\},$$
the second condition because of determinant 1. From the multiplicativity of
the norm, the second relation can be written
$N_{D|K}(a)=N_{D|K}(\overline{a})$.  Since the norm is given by the
determinant, this means, in $D^*({\Bbb R})$,
$$L({\Bbb R})\cong \left\{ a\in D^*({\Bbb R})\Big| a\not\in K\Ra \det(a)\in
  {\Bbb R}^*\right\} \cong U(1)\cdot GL^0(d,{\Bbb C}),$$ where the $U(1)$
factor arises from the elements in $U(1)\cap K$. Similarly, the unipotent
radical takes the form
$$U \cong \left\{\left(\begin{array}{cc} 1 & b \\ 0 & 1 \end{array}\right)
  \Big| b\in D,\ b+\overline{b}=0\right\} \cong \sqrt{-\eta}D^+,$$ which is
$d^2$-dimensional. An explicit isomorphism $U({\Bbb R})\cong {\cal
  H}_d({\Bbb C})$ is given by the identity on $D^+({\Bbb R})$, since $x\in
D^+\Ra x=\overline{x}$, hence as a real matrix, this element is hermitian.
\end{itemize}
The maximal compact subgroups of $G_D({\Bbb R})$ are the groups ${\cal
  K}({\Bbb R})$, where ${\cal K}$ is the group of (\ref{eCom}). Indeed, it
is easy to see that the stabiliser of the non-isotropic vector $(1,1)$ is
the stabiliser, in the domain $\cD$, of the base point, which is by
definition the maximal compact subgroup.

Now consider the subgroups $U(L^2,h)\inn G_D$; let us introduce the
notation $G_L$ for these subgroups. Then, since $G_L$ is of the type $d=1$
in Proposition \ref{p9.1}, it follows that Theorem \ref{t9.1} applies to
$G_L({\Bbb R})$ as well as to $G_D({\Bbb R})$.
\begin{proposition}\label{p10.1} We have a commutative diagram
  $$\begin{array}{ccc} G_L({\Bbb R}) & \hra & G_D({\Bbb R}) \\ \downarrow &
    & \downarrow \\ \cD_L & \hra & \cD_D, \end{array}$$ where $\cD_L$ and
  $\cD_D$ are hermitian symmetric domains, and the maps $G({\Bbb R})\lra
  \cD$ are the natural projections. Moreover, the subdomains $\cD_L$ are as
  follows:
\begin{itemize}\item[1)] $d=1$; We now assume that the extension $k/k'$ is
  {\it Galois}. Then, if $\deg_{{\Bbb Q}}k'=f'$, $f/f'=m$, we have
  $\cD_L\cong (\HH)^{f'}$ and $\cD_D\cong ((\HH)^m)^{f'}$ and the embedding
  $\cD_L\inn \cD_D$ is given by $\HH\hra (\HH)^m$ diagonally, and the
  product of this $f'$ times.
\item[2)] $d=2$; $\cD_L\cong \left(\begin{array}{cc}\tau_1 & 0 \\ 0 &
      b^{\zeta_1}\tau_1
\end{array}\right)\times \cdots \times \left(\begin{array}{cc}\tau_1 & 0 \\ 0 &
  b^{\zeta_f}\tau_1
\end{array}\right)$, where $\zeta_i:k\lra {\Bbb R}$ denote the distinct real
embeddings of $k$.
\item[3)] $d\geq 3$; $\cD_L\cong \left(\begin{array}{ccc}\tau_1 & & 0 \\ &
      \ddots & \\ 0 & & \tau_d\end{array}\right)^f$.
\end{itemize}
\end{proposition}
{\bf Proof:} Considering the tube realisation of the domain $\cD_D$, we
have the base point
$$o=(\diag(i,i,\ldots,i))^f,$$ and we just calculate the orbit of the
${\Bbb R}$-group $G_L({\Bbb R})$ of this point; this gives the
corresponding subdomain. The commutativity of the diagram follows from
this.
\begin{itemize}\item[1)] $d=1$: Since $k/k'$ is Galois, we have that
  $G_L\inn G_D$ is the subgroup fixed under the natural $Gal_{k|k'}$
  action; hence lifting the ${\Bbb Q}$-groups to ${\Bbb R}$, we get the
  identification
  $$(Res_{k'|{\Bbb Q}}G_L)_{{\Bbb R}} \cong (Res_{k|{\Bbb Q}}G_D)_{{\Bbb
      R}}^{Gal_{k|k'}}.$$ In terms of the domains, if $\cD_D=(\HH)^f$, then
  $Gal_{k|k'}$ acts by permuting $m$ copies at a time, and the invariant
  part is the diagonally embedded $\HH\hra (\HH)^m,\ z\mapsto
  (z,\ldots,z)$. Since $\cD_L\cong (\HH)^{f'}$, the result follows in this
  case.
\item[2)] $d=2$: First we remark that $G_{{\Bbb R}}$ is conjugate to the
  {\it standard} $Sp(4,{\Bbb R})$ (by which we mean the symplectic group
  with respect to the symplectic form $J={0\ \ 1 \choose -1\ 0}$) under the
  element
  $$q=\left(\begin{array}{c|c}\begin{array}{cc} 1 & 0 \\ 0 & 1 \end{array}
      & 0 \\ \hline 0 & \begin{array}{cc} 0 & 1 \\ -1 & 0
\end{array}\end{array}\right),$$
that is $qG_{{\Bbb R}}q^{-1}=Sp(4,{\Bbb R})$. To see this, we recall how
the ($D$-valued) hyperbolic (hermitian) form is related to alternating
forms.  Since $D_{{\Bbb R}}\cong M_2({\Bbb R})$, our vector space $D^2$
over ${\Bbb R}$ is $V_{{\Bbb R}}=(M_2({\Bbb R}))^2$; the hyperbolic form
gives a $M_2({\Bbb R})$-valued form on $V_{{\Bbb R}}$. We have matrix units
$e_{ij},\ i,j=1,2$ in $M_2({\Bbb R})$, and setting $W_1=e_{11}V,\
W_2=e_{22}V$, the spaces $W_i$ are two-dimensional vector spaces over
${\Bbb R}$. The hermitian form $h({\bf x},{\bf y})$,
restricted to $W_1$, turns out to
have values in ${\Bbb R} e_{12}$: if
$${\bf x}=\left(\left(\begin{array}{cc} x_{11} & x_{12} \\ x_{21} &
      x_{22}\end{array}\right),\left(\begin{array}{cc} x'_{11} & x'_{12} \\
      x'_{21} & x'_{22}\end{array}\right)\right), \ \ \ {\bf
  y}=\left(\left(\begin{array}{cc} y_{11} & y_{12} \\ y_{21} &
      y_{22}\end{array}\right),\left(\begin{array}{cc} y'_{11} & y'_{12} \\
      y'_{21} & y'_{22}\end{array}\right)\right),$$ then ${\bf x, y}\in W_1
\iff x_{2i}=x'_{2i}=y_{2i}=y'_{2i}=0$. Then we have\footnote{the involution
  on $M_2(k)$ is as in (\ref{e2.2}).}
$$h({\bf x},{\bf y})=\left(\begin{array}{cc} x_{11} & x_{12} \\ 0 & 0
\end{array}\right)\left(\begin{array}{cc} 0 & -y'_{12} \\ 0 & y'_{11}
\end{array} \right) + \left(\begin{array}{cc} x'_{11} & x'_{12} \\ 0 & 0
\end{array}\right) \left( \begin{array}{cc} 0 & -y_{12} \\ 0 & y_{11}
\end{array} \right) $$
$$ = \left(\begin{array}{cc} 0 & -x_{11}y'_{12} + x_{12}y'_{11}
    -x'_{11}y_{12} + x'_{12}y_{11} \\ 0 & 0 \end{array}\right).$$ Then the
alternating form $\langle {\bf x},{\bf y} \rangle$ is defined by
$$h({\bf x},{\bf y})=\langle {\bf x}, {\bf y}\rangle\cdot e_{12},\ {\bf
  x,y}\in W_1.$$ In other words, viewing the symplectic group as a subgroup
of $GL(2,M_2({\Bbb R}))$, the symplectic form will be given by $J_1={0\ \
  J_2 \choose J_3\ 0}$, where $J_2$ and $J_3$ are symplectic forms on
${\Bbb R}^2$.  It is clear that $q{0\ \ 1_2 \choose 1_2\ 0}{^tq}=J_1$ is of
this form, and this shows $qG_{{\Bbb R}}q^{-1}=Sp(4,{\Bbb R})$.

Next, recalling that we view $D$ as a cyclic algebra of $2\times 2$
matrices, an element $A={\alpha\ \beta \choose \gamma\ \delta}\in M_2(D)$
is actually a $4\times 4$ matrix; each of $\alpha,\beta,\gamma, \delta$ are
of the form (\ref{e2.1}), and this $4\times 4$ matrix must be conjugated by
$q$, then the element acts in the well known manner on the domain $\cD_D$:
$$M=\left(\begin{array}{cc} A & B \\ C & D \end{array}\right) \in
Sp(4,{\Bbb R}),\ \ \ \tau=\left(\begin{array}{cc} \tau_1 & \tau_{12} \\
    \tau_{12} & \tau_2\end{array}\right)\in \cD_D,\ \ \ M\cdot \tau =
(A\tau + B)(C\tau +D)^{-1}.$$ Finally we note that the subgroup
$U(L^2,h)\inn U(D^2,h)$ consists of the elements of the form as in
Proposition \ref{p7.1} 2):
$$x=\left(\begin{array}{cc} \alpha & 2\beta/\sqrt{-ab} \\
    \gamma\sqrt{-ab}/2 & \delta
\end{array}\right),\ \ \ \alpha,\beta,\gamma,\delta\in k,$$
and the element $\sqrt{-ab}=ec$ as in (\ref{e7.1}), so that
$(\sqrt{-ab})^{-1}=(ec)^{-1}$. In other words, as a $4\times 4$ matrix, the
element $x$ above is
\begin{equation}\label{e10b.1} x=\left(\begin{array}{cc|cc}
      \alpha & 0 & 0 & {2\beta \over b\sqrt{a}} \\ 0 & \alpha & {2\beta
        \over -\sqrt{a}} & 0 \\ \hline 0 & {-\gamma\sqrt{a}\over 2} &
      \delta & 0 \\ {b\gamma\sqrt{a}\over 2} & 0 & 0 & \delta \\
\end{array}\right),
\end{equation}
and after conjugating with $q$ this becomes \begin{equation}\label{e10b.2}
  qxq^{-1}=\left(\begin{array}{cc|cc} \alpha & 0 & {2\beta \over b\sqrt{a}}
      & 0 \\ 0 & \alpha & 0 & {2\beta\over \sqrt{a}} \\ \hline
      {b\gamma\sqrt{a}\over 2} & 0 & \delta & 0 \\ 0 & {\gamma\sqrt{a}\over
        2} & 0 & \delta\end{array}\right).
\end{equation}
{}From the form of the matrix (\ref{e10b.2}) we see that the corresponding
subdomain is contained in the set of diagonal matrices ${\tau_1\ \ 0\choose
  0\ \tau_2}^f$. Now calculating, setting $g=qxq^{-1}$, we have
\begin{equation}\label{e10b.2a} g\left(\begin{array}{cc}\tau_1 & 0 \\ 0 &
      \tau_2\end{array}\right)^f=\left(\left(\begin{array}{cc}
        \left({\alpha b \tau_1 +{2\beta\sqrt{a}\over a} \over
            {\gamma\sqrt{a}\over 2}b\tau_1+\delta}\right){1\over b} & 0 \\
        0 & {\alpha \tau_2 +{2\beta\sqrt{a}\over a} \over
          {\gamma\sqrt{a}\over 2}\tau_2+\delta}
      \end{array}\right)^{\zeta_1},\ldots, \left(\begin{array}{cc}
      \left({\alpha
            b \tau_1 +{2\beta\sqrt{a}\over a} \over {\gamma\sqrt{a}\over
              2}b\tau_1+\delta}\right){1\over b} & 0 \\ 0 & {\alpha \tau_2
          +{2\beta\sqrt{a}\over a} \over {\gamma\sqrt{a}\over
            2}\tau_2+\delta}
      \end{array}\right)^{\zeta_f} \right),
\end{equation}
So in particular, $G_L$ maps the set of matrices
\begin{equation}\label{e10b.2b} \cD_L=\left(\begin{array}{cc}\tau_1 & 0 \\
      0 & b^{\zeta_1}\tau_1 \end{array}\right)\times \cdots \times
  \left(\begin{array}{cc}\tau_1 & 0 \\ 0 & b^{\zeta_f}\tau_1
    \end{array}\right)
\end{equation}
into itself. (In the expression (\ref{e10b.2a}), the bracket
$(...)^{\zeta_i}$ means $\zeta_i$ is applied to all matrix elements).  This
proves the theorem in the case $d=2$.
\begin{remark} The description of the subdomain $\cD$ in terms of $b$,
  depends on the way we described $D$ as a cyclic algebra. We took
  $\ell=k(\sqrt{a}),\ D=(\ell/k,\sigma,b)$. But we could also consider
  $\ell'=k(\sqrt{b}),\ D'=(\ell'/k,\sigma',a)$; these two cyclic algebras
  are isomorphic, an isomorphism being given by extending the identity on
  the center $k$ by
\begin{eqnarray*} D & \stackrel{\sim}{\lra} & D' \\ e & \mapsto & c' \\
  c & \mapsto & e',\end{eqnarray*} where $$e=\left(\begin{array}{cc} 0 & 1
    \\ b & 0 \end{array}\right),\quad c=\left(\begin{array}{cc} \sqrt{a} &
    0 \\ 0 & -\sqrt{a} \end{array}\right), \hspace*{2cm}
e'=\left(\begin{array}{cc} 0 & 1 \\ a & 0\end{array}\right), \quad
c'=\left( \begin{array}{cc} \sqrt{b} & 0 \\ 0 &
    -\sqrt{b}\end{array}\right).$$ Indeed, as one sees immediately,
$(ec)^2=-ab=(e'c')^2$, so the map above is a morphism, which is clearly
bijective.  Consequently, it would be better to denote (\ref{e10b.2b})
$\cD_{L,b}$, since we {\it also} have a subdomain
$$\cD_{L,a}=\left(\begin{array}{cc} \tau_1 & 0 \\ 0 &
    a^{\zeta_1}\tau_1\end{array}\right)\times \cdots \times
\left(\begin{array}{cc} \tau_1 & 0 \\ 0 &
    a^{\zeta_f}\tau_1\end{array}\right),$$ neither of which is a priori
privledged.
\end{remark}

\item[3)] $d\geq 3$: Here the situation is simpler; the matrix $H$ is the
  matrix of the hermitian form whose symmetry group acts on the usual
  unbounded realisation of the domain of type $\bf I_{d,d}$, hence no
  conjugation is necessary. Let $G_L$ be the subgroup we wish to consider,
  and recall from Proposition \ref{p7.1} that it consists of matrices
  $$x=\left(\begin{array}{cc} \alpha & 2\beta/\sqrt{-\eta} \\
      \gamma\sqrt{-\eta}/2 & \delta
\end{array}\right),\ \ \ \alpha,\beta,\gamma,\delta\in k,\
  \alpha\delta-\gamma\beta=1.$$
Thinking of $D$ itself as $d\times d$ matrices, this element is
\begin{equation}x=\left(\begin{array}{cccc|cccc} \alpha & & & 0 & {2\beta \over
        \sqrt{-\eta}} & & & 0\\ & \alpha^{\sigma} & & & & {2\beta^{\sigma}
        \over \sqrt{-\eta}} & & \\ & & \ddots & & & & \ddots & \\ 0 & & &
      \alpha^{\sigma^{d-1}} & 0 & & & {2\beta^{\sigma^{d-1}} \over
        \sqrt{-\eta}} \\ \hline {\gamma\sqrt{-\eta} \over 2} & & & 0 &
      \delta & & & 0 \\ & {\gamma^{\sigma}\sqrt{-\eta} \over 2} & & & &
      \delta^{\sigma} & & \\ & & \ddots & & & & \ddots & \\ 0 & & &
      {\gamma^{\sigma^{d-1}}\sqrt{-\eta} \over 2} & 0 & & &
      \delta^{\sigma^{d-1}} \end{array} \right),
\end{equation}
and the orbit of the base point $\diag(i,\ldots, i)^f$ under elements as
$x$ is clearly the diagonal subdomain $(\diag(\tau_1,\ldots, \tau_d))^f$,
and the action is, for $k={\Bbb Q}$,
\begin{equation}\label{e10b.5} \diag(\tau_1,\ldots,\tau_d) \mapsto \left(
    {\alpha\tau_1+{2\beta \over \sqrt{-\eta}} \over
      {\gamma\sqrt{-\eta}\over 2} \tau_1 + \delta},\ldots ,
    {\alpha^{\sigma^{d-1}}\tau_{d}+{2\beta^{\sigma^{d-1}}\over
        \sqrt{-\eta}} \over {\gamma^{\sigma^{d-1}}\sqrt{-\eta}\over
        2}\tau_d + \delta^{\sigma^{d-1}} }\right)=:\left(\begin{array}{cc}
      \alpha & {2\beta\over \sqrt{-\eta}} \\ {\gamma\sqrt{-\eta}\over 2} &
      \delta\end{array} \right)(\tau_1,\ldots,\tau_d),
\end{equation}
and for general $k$, letting as above $\zeta_1,\ldots,\zeta_f$ denote the
$f$ embeddings of $k$,
\begin{equation}\label{e10b.6} (\diag(\tau_1,\ldots,\tau_d))^f \mapsto
  \left( \left(\begin{array}{cc} \alpha & {2\beta\over \sqrt{-\eta}} \\
        {\gamma\sqrt{-\eta}\over 2} & \delta\end{array}
    \right)^{\zeta_1}(\tau_1,\ldots,\tau_d),\ldots, \left(\begin{array}{cc}
        \alpha & {2\beta\over \sqrt{-\eta}} \\ {\gamma\sqrt{-\eta}\over 2}
        & \delta\end{array} \right)^{\zeta_f}(\tau_1,\ldots,\tau_d)
  \right).
\end{equation}
\end{itemize}
This establishes all statements of the theorem. \ende

Finally, let $F$ be the rational boundary component $\in \cD_D^*$
corresponding to the isotropic vector $(0,1)$, of which the parabolic $P$
of (\ref{e6.1}) is the normaliser, $P({\Bbb R})=N(F)$. Note that by
construction, $(0,1)$ is also isotropic for $G_L$; let $F_L$ be the
rational boundary component $\in \cD_L^*$ corresponding to it, so the
parabolic in $G_L$, $P_L$ of (\ref{e6.1}) is the normaliser, $P_L({\Bbb
  R})=N(F_L)$. Then clearly
\begin{equation}\label{e10.2} P_L=G_L\cap P.
\end{equation}
Also we have
\begin{proposition}\label{p10.2} Let $\cD_L\inn \cD_D,\ F_L\in \cD_L^*,\
  F\in \cD_D^*$ be as above, and let $i:\cD_L^*\hra \cD_D^*$ be the induced
  inclusion of Satake compactifications of the domains. Then $i(F_L)=F$.
\end{proposition}
We will refer to this as ``the subdomain $\cD_L$ contains $F$ as a rational
boundary component.'' The proposition itself is, after unraveling the
definitions, nothing by (\ref{e10.2}).

\section{Arithmetic groups}
It is well-known how to construct arithmetic groups $\Gamma\inn G_D(K)$; we
recall this briefly. View $D^2$ as a $D$-vector space, and let $\Delta\inn
D$ be a maximal order, that is an ${\cal O}_K$-module (where $D$ is central
simple over $K$, $K=k$ if $d=2$), which spans $D$ as a $K$-vector space, is
a sub{\it ring} of $D$, and is maximal with these properties. Then set
\begin{equation} \label{e11.1}\Gamma_{\Delta}:=G_D(K)\cap M_2(\Delta).
\end{equation}
If we view $\Delta^2\inn D^2$ as a lattice, then this is also described as
$$\Gamma_{\Delta}=\left\{ g\in G_D(K) \Big| g(\Delta^2)\inn \Delta^2
\right\}.$$ {}From the second description it is clear that
$\Gamma_{\Delta}$ is an {\it arithmetic} subgroup. We will denote the
quotient of the domain $\cD_D$ by this arithmetic subgroup by
\begin{equation}\label{e11.2} X_{\Gamma_{\Delta}}=X_{\Delta}=
  \Gamma_{\Delta}\bs \cD_D.
\end{equation}
If $\Gamma_{\Delta}$ were fixed-point free on $\cD_D$ this would be a
non-compact complex manifold; if $\Gamma_{\Delta}$ has fixed points (as it
usually does), these yield singularities of the space $X_{\Delta}$, which
is still an analytic variety. It is a $V$-variety in the language of
Satake. In fact it has a {\it global} finite smooth cover: let $\Gamma\inn
\Gamma_{\Delta}$ be a normal subgroup of finite index which does act freely
(such exist by a theorem of Selberg). Then we have a Galois cover
$$X_{\Gamma}\lra X_{\Gamma_{\Delta}},$$ and $X_{\Gamma}$ is smooth. Hence
the singularities of $X_{\Gamma_{\Delta}}$ are encoded in the action of the
finite Galois group of the cover.

The Baily-Borel compactification of $X_{\Gamma}$ ($\Gamma\inn
\Gamma_{\Delta}$) is an embedding
$$X_{\Gamma}\inn X_{\Gamma}^*\inn {\Bbb P}^N,$$ where $ X_{\Gamma}^*$ is a
normal algebraic variety. Since the boundary components of $\cD_D$ are
points (Theorem \ref{t9.1}), it follows that the singularities of $
X_{\Gamma}^*$ contained in the boundary $ X_{\Gamma}^*-X_{\Gamma}$ are also
{\it isolated points}. In particular, if $\Gamma\inn \Gamma_{\Delta}$ has
no elements of finite order, then $ X_{\Gamma}^*$ is smooth outside a
finite set of isolated points, all contained in the compactification locus.

The singularities of $ X_{\Gamma}^*$ may be resolved by means of toroidal
compactifications; let us denote such a compactification by
$\overline{X}_{\Gamma}$, for which we make the assumptions:
\begin{itemize}\item[i)] $\overline{X}_{\Gamma}-X_{\Gamma}$
  is a normal crossings divisor,
\item[ii)] $\overline{X}_{\Gamma}\lra X_{\Gamma}^*$ is a resolution of
  singularities and $\overline{X}_{\Gamma}$ is projective algebraic.
\end{itemize}
Such a compactification $\overline{X}_{\Gamma}$ depends on a set of
polyhedral cones, and one can make choices such that i) and ii) are
satisfied, by the results of \cite{SC}. We will not need these in detail;
the mere existence will be sufficient.

Recall that a cusp of an arithmetic group $\Gamma$ is a maximal flag of
boundary components of $\Gamma$; in the case at hand this is just a point
in $ X_{\Gamma}^*-X_{\Gamma}$. The {\it number of cusps} is the number of
points in $ X_{\Gamma}^*- X_{\Gamma}$, and may be defined in the following
ways:
\begin{itemize}\item[a)] It is the number of $\Gamma$-equivalence classes of
  parabolic subgroups of $G_D$ or $SG_D$; here $\Gamma$ acts on the
  parabolic subgroups by conjugation.
\item[b)] It is the number of $\Gamma$-equivalence classes of isotropic
  vectors of the hermitian plane $(D^2,h)$; here $\Gamma$ is acting as a
  subgroup of $G_D$ on the vector space.
\end{itemize}
In the $d=1$ case we are talking about the number of cusps of a Hilbert
modular variety. Let us recall how these numbers are determined. For this
we consider first the group $SL_2(k)$ (instead of $SU(K^2,h)$), which acts
on a product of upper half planes, $\cD_D=\HH^f$. The boundary is the
product $({\Bbb P}^1({\Bbb R}))^f$, and the rational boundary components
are the points of the image ${\Bbb P}^1(k)\hra ({\Bbb P}^1({\Bbb R}))^f$
given by $x\mapsto (x^{\zeta_1},\ldots,x^{\zeta_f})$. Hence we may denote
the boundary components by $\xi=(\xi_1:\xi_2),\ \xi_i\in {\cal O}_k$. Let
$\aa_{\xi}$ denote the ideal generated by $\xi_1,\xi_2$. Then one has
\begin{proposition}\label{p13.1} Two boundary components
  $\xi=(\xi_1:\xi_2)$ and $\eta=(\eta_1:\eta_2)$ are equivalent under
  $SL_2({\cal O}_k)$ $\iff$ the ideals $\aa_{\xi}$ and $\aa_{\eta}$ are in
  the same class. In particular, the number of cusps is the class number of
  $k$.
\end{proposition}
{\bf Proof:} One direction is trivial: if ${a\ b\choose c\ d}{\xi_1\choose
  \xi_2}={\eta_1\choose \eta_2},\ g={a\ b\choose c\ d}\in SL_2({\cal
  O}_k)$, then the ideals $\aa_{\xi}$ and $\aa_{\eta}$ are the {\it same}.
This is because $g$ is unimodular, so affects just a change of base in the
ideal $\aa_{\xi}$. Conversely, suppose $\xi$ and $\eta$ are given, and
suppose that $\aa_{\xi}$ and $\aa_{\eta}$ have the same class. After
multiplication by an element $c\in k^*$, we may assume
$\aa_{\xi}=\aa_{\eta}=\aa$.  Furthermore, writing ${\cal
  O}_k=\aa\aa^{-1}=\xi_1\aa^{-1}+\xi_2\aa^{-1}$, we see that $1\in {\cal
  O}_k$ can be written
\begin{equation}\label{e13.1}
  1=\xi_1\xi_2'-\xi_2\xi_1'=\eta_1\eta_2'-\eta_2\eta_1',\ \ \xi_i',
  \eta_i'\in \aa^{-1}.
\end{equation}
But this means the matrices (acting from the left)
\begin{equation}\label{e13.2} M_{\xi}=\left(\begin{array}{cc} \xi_1 & \xi_1'
      \\ \xi_2 & \xi_2' \end{array}\right),\quad\quad
  M_{\eta}=\left(\begin{array}{cc} \eta_1 & \eta_1' \\ \eta_2 & \eta_2'
\end{array}\right),
\end{equation}
are in $SL_2(k)$, and
\begin{equation}\label{e13.3} M_{\xi}{ 1 \choose
    0}={\xi_1\choose \xi_2}, \quad\quad M_{\eta}{1\choose 0}={\eta_1\choose
    \eta_2}.
\end{equation}
Of course $M_{\xi},\ M_{\eta}$ are not in $SL_2({\cal O}_k)$, but from the
fact that $\xi_i',\eta_i'\in \aa^{-1}$ it follows that
$M_{\xi}M_{\eta}^{-1}\in SL_2({\cal O}_k)$. Hence
$$(M_{\xi}M_{\eta}^{-1}){\eta_1\choose \eta_2} ={\xi_1\choose \xi_2},$$ and
the cusps are conjugate under $SL_2({\cal O}_k)$. \ende We now translate
this into the corresponding statement for the hyperbolic plane $(K^2,h)$,
where $K|k$ is an imaginary quadratic extension. First note that if
${\xi_1\choose \xi_2}\in K^2$ is {\it isotropic}, then
$$Tr_{K|k}(\xi_1\overline{\xi}_2)=
 \xi_1\overline{\xi}_2+\xi_2\overline{\xi}_1=0.$$
In fact, the converse also holds,
\begin{lemma}\label{l14.1} A vector ${\xi_1\choose \xi_2}\in K^2$ is
  isotropic with respect to $h$ $\iff
  Tr(\xi_1\overline{\xi}_2)=Tr(\overline{\xi}_1\xi_2)=0$ $\iff$ $\xi_2=0$
  or $Tr(\xi_1\xi_2^{-1})=0$.
\end{lemma}
{\bf Proof:} By definition, ${\xi_1\choose \xi_2}$ is isotropic $\iff$
$\xi_1\overline{\xi}_2+\xi_2\overline{\xi}_1=0$, but this is
$Tr(\xi_1\overline{\xi}_2)=0$.  Noting that $\xi_2^{-1}={1\over
  N(\xi_2)}\overline{\xi}_2$, if $\xi_2\neq0$, we get the second
equivalence, from the linearity of the trace,
$Tr(\xi_1\xi_2^{-1})=Tr(\xi_1{1\over N(\xi_2)}\overline{\xi}_2)= {1\over
  N(\xi_2)}Tr(\xi_1\overline{\xi}_2)$. If $\xi_2=0$, then
${\xi_1\choose\xi_2}={\xi_1\choose 0}$ is isotropic anyway, establishing
both equivalences as stated. \ende Let $K^0$ denote the purely imaginary
(traceless) elements of $K$, $K^0\cong \sqrt{-\eta}k$ for
$K=k(\sqrt{-\eta})$. It follows from Lemma \ref{l14.1} that any $x\in K^0$
of the form $x=\xi_1\xi_2^{-1}$ or $x=\xi_1\overline{\xi}_2$ determines an
isotropic vector ${\xi_1\choose \xi_2}$ of $K^2$. Consider in particular
the case that ${\xi_1\choose \xi_2}\in {\cal O}_K^2$ is integral.  Recall
that (non-zero) $\xi_1$ and $\xi_2$ are relatively prime $\iff$ there are
$x,y\in {\cal O}_K$ such that $x\xi_1+y\xi_2=1$. Define the set of
relatively prime integral isotropic vectors,
\begin{equation}\label{e14.2} \scI=\left\{ {\xi_1\choose \xi_2}\in
    {\cal O}_K^2-(0,0) \Big| Tr(\xi_1\overline{\xi}_2)=0; \hbox{ and if }
    \xi_i\neq 0, i=1,2,\ \hbox{ then } \exists_{x,y\in {\cal O}_K}\
    x\xi_1+y\xi_2=1\right\}.
\end{equation}
Following the proof of Proposition \ref{p13.1}, we consider when two
isotropic integral vectors are conjugate under $S\Gamma_{\Delta}=SU({\cal
  O}_K^2,h)$.  Let now $\aa_{\xi}$ be the ideal generated by $\xi_1,\xi_2$
in ${\cal O}_K$.
\begin{proposition}\label{p14.1} Two isotropic vectors $\xi={\xi_1\choose
    \xi_2}$ and $\xi'={\xi'_1\choose \xi'_2}\in \scI$ are conjugate under
  $SU({\cal O}_K^2,h)$ $\iff$ the ideals $\aa_{\xi}$ and $\aa_{\xi'}$ are
  equivalent in $K$ $\iff$ the ideals $N(\aa_{\xi})$ and $N(\aa_{\xi'}) $
  are equivalent in $k$.
\end{proposition}
{\bf Proof:} If $g{\xi_1\choose \xi_2}={\xi_1'\choose \xi_2'},\ g\in
SU({\cal O}_K^2,h)$, then since $g$ is unimodular, the ideal classes
coincide.  Conversely, suppose ${\xi_1\choose \xi_2}$ and
${\xi_1'\choose{\xi'_2}}$ are equivalent; again we may assume
$\aa_{\xi}=\aa_{\xi'}=\aa$. Then, as above, there are $\rho,\ \rho'\in
(\aa^{-1})^2$ such that (\ref{e13.1}) holds.  Consequently we have
$M_{\xi}$, $M_{\xi'}$, but now in $SU(K^2,h)$, such that (\ref{e13.2}) and
(\ref{e13.3}) hold. It follows that $M_{\xi}M_{\xi'}^{-1}\in SU({\cal
  O}_K^2,h)$ maps $\xi'$ to $\xi$. This completes verification of the first
$\iff$.  The second equivalence then follows from Proposition \ref{p5.1}
and Proposition \ref{p13.1}. \ende We now derive an analogue of the above
for general $D$ as in (\ref{e1.0}).  We first note that Lemma \ref{l14.1}
is true here also.
\begin{lemma}\label{l15.1} Let $\xi=(\xi_1,\xi_2)\in D^2$ be given. Then
  $\xi$ is isotropic $\iff$
  $\xi_1\overline{\xi}_2+\xi_2\overline{\xi}_1=0$.
\end{lemma}
{\bf Proof:} We have
$h(\xi,\xi)=\xi_1\overline{\xi}_2+\xi_2\overline{\xi}_1=0$, so that
$h(\xi,\xi)=0 $ if and only if
$\xi_1\overline{\xi}_2+\xi_2\overline{\xi}_1=0$.  \ende Next we note that
any isotropic vector is conjugate in $G_D$ to the standard isotropic vector
$(0,1)$.
\begin{lemma}\label{l15.2} Let $\xi=(\xi_1,\xi_2)\in D^2$ be isotropic.
  Then there is a matrix $M_{\xi}\in G_D$ such that $(0,1)M_{\xi}=\xi$.
\end{lemma}
{\bf Proof:} If $M_{\xi}=\left(\begin{array}{cc} \xi_1' & \xi_2' \\ \xi_1
    &\xi_2\end{array}\right)$, then $(0,1)M_{\xi}=\xi$. So we must show the
existence of such an $M_{\xi}\in G_D$. The equations to be solved (for
$\xi_i'$) are
\begin{equation}\label{e15.1}\begin{minipage}{15cm}\begin{itemize}\item[1)]
      $\xi_1'\overline{\xi}_2+\xi_2'\overline{\xi}_1=1$.
    \item[2)] $\xi_1'\overline{\xi}'_2+\xi_2'\overline{\xi}'_1=0.$
    \item[3)] $\xi_1\overline{\xi}_2+\xi_2\overline{\xi}_1=0.$
\end{itemize}
\end{minipage}\end{equation}
Since $\xi$ is isotropic, 3) is fulfilled by Lemma \ref{l15.1}. If
$\xi_1\neq 0$, then setting $\xi_1'=0, \xi_2'=\overline{\xi}_1^{-1}$,
$\xi'$ fulfills both 1) and 2). If $\xi_1=0$, then $\xi_2\neq 0$ and we set
similarly $\xi_1'=\overline{\xi}_2^{-1}$ and $\xi_2'=0$. \ende Now we
consider as above integral isotropic vectors. Let $\Delta\inn D$ be a
maximal order, and set
$$\scI_{\Delta}=\left\{ (\xi_1,\xi_2)\in \Delta^2-\{(0,0)\} \Big|
  \xi_1\overline{\xi}_2+\xi_2\overline{\xi}_1=0, \xi_1, \xi_2\neq 0 \Ra
  \exists_{x,y\in \Delta}\ \xi_1x+\xi_2y=1\right\}.$$
\begin{lemma} \label{l15.3} Let $(\xi_1,\xi_2)\in \scI_{\Delta}$. Then there
  exists a matrix
  $$M_{\xi}=\left(\begin{array}{cc} \xi_1'' & \xi_2'' \\ \xi_1
      &\xi_2\end{array}\right)\in \Gamma_{\Delta}.$$
\end{lemma}
{\bf Proof:} If $\xi_2=0$, then $\xi''_2=1$, any $\xi''_1$ with
$\xi''_1+\overline{\xi}''_1=0$ gives such a matrix, similarly, if $\xi_1=0$,
then
$\xi''_1=1$ gives a solution. So we assume $\xi_i\neq 0$; then by assumption
we have $(x,y)\in \Delta^2$ such that $\xi_1x+\xi_2 y=1$, hence
$\overline{x}\overline{\xi}_1+\overline{y}\overline{\xi}_2=1$. So setting
$\xi_2'=\overline{x}, \xi_1'=\overline{y}$, the first relation of
(\ref{e15.1}) is satisfied. Again 3) is satisfied by assumption, so we must
verify 2). It turns out that we may have to alter $\xi_i'$ to achieve this.
The relation 1) can be expressed as $h(\xi',\xi)=1$, where
$\xi=(\xi_1,\xi_2),\ \xi'=(\xi_1',\xi_2')$. Since $\xi$ is isotropic by
Lemma \ref{l15.1}, $h(\xi,\xi)=0$, and with respect to the base
$<\xi,\xi'>$ of $D^2$, $h$ is given by a matrix $H_{\xi,\xi'}={0\ 1 \choose
  1\ \ge}$, where $\ge=h(\xi',\xi')= \delta + \overline{\delta}$,
$\delta=\xi_1'\overline{\xi}'_2$ (if $\ge\neq0$; otherwise we are done), in
particular $\delta\in \Delta$. Now setting
$$\xi''=(\xi_1'',\xi_2'')=(-\xi_1\overline{\delta}+\xi_1',
 -\xi_2\overline{\delta}+\xi_2')\in
\Delta^2,$$ we can easily verify $h(\xi'',\xi'')=0,\ \
h(\xi,\xi'')=h(\xi'',\xi)=1$, so this vector $\xi''$ gives a matrix
$M_{\xi} \in \Gamma_{\Delta}$, as was to be shown. \ende We require the
following refinement of Lemma \ref{l15.2}. For this, given
$(\xi_1,\xi_2)\in \Delta^2$, let $\aa_{\xi}$ denote the left ideal in
$\Delta$ generated by the elements $\xi_1, \xi_2$.
\begin{lemma}\label{l16.1} Given $\xi=(\xi_1,\xi_2)\in \Delta^2$ isotropic,
  the entries $\xi_i'$ of the matrix $M_{\xi}$ of Lemma \ref{l15.2} are in
  $\aa_{\xi}^{-1}$.
\end{lemma}
{\bf Proof:} This follows from the proof of \ref{l15.2}, as we took just
inverses of elements $\xi_i$. \ende We can now consider when two integral
isotropic vectors are equivalent under $\Gamma_{\Delta}$.
\begin{proposition}\label{p16.1} Let $\xi=(\xi_1,\xi_2),\
  \eta=(\eta_1,\eta_2)\in \Delta^2$ be two integral isotropic vectors,
  $\aa_{\xi},\ \aa_{\eta}$ the left ideals generated by the components of
  $\xi$ and $\eta$, respectively. Then $\xi$ and $\eta$ are equivalent
  under an element of $\Gamma_{\Delta}$ $\iff$ $\aa_{\eta}\cong \aa_{\xi}$
  as $\Delta$-modules.
\end{proposition}
{\bf Proof:} If there is an element $g={a\ b\choose c\ d}\in
\Gamma_{\Delta}$ with $(\xi_1,\xi_2)g=(\eta_1,\eta_2)$, then
$\eta_1=\xi_1a+\xi_2c,\ \eta_2=\xi_1b+\xi_2d$, hence $\aa_{\eta}\inn
\aa_{\xi}$ and the map
\begin{eqnarray*} \phi:\aa_{\eta} & \lra & \aa_{\xi} \\
  \eta_1\lambda_1+\eta_2\lambda_2 & \mapsto & \xi_1(a\lambda_1+b\lambda_2)+
  \xi_2(c\lambda_1+d\lambda_2)
\end{eqnarray*}
is an isomorphism of $\Delta$-modules. The inverse is given is a similar
manner by $g^{-1}$. Conversely, if $\aa_{\eta}\cong \aa_{\xi}$ we may
assume $\aa_{\xi}=\aa_{\eta}$, hence also $\aa_{\xi}^{-1}=
\aa_{\eta}^{-1}$.  Applying Lemma \ref{l16.1} and making use of
$\aa_{\xi}\aa_{\xi}^{-1}=\Delta$, which holds in central simple division
algebras over number fields, we find for the matrices $M_{\xi},\ M_{\eta}$
associated by Lemma \ref{l15.2} to $\xi$ and $\eta$, respectively, that
$M_{\xi}M_{\eta}^{-1}\in \Gamma_{\Delta}.$ For example, if $\xi_1\neq 0,\
\eta_1\neq 0$, then
$$M_{\xi}M_{\eta}^{-1} = \left(\begin{array}{cc}0 & \overline{\xi}_1^{-1}
    \\ \xi_1 & \xi_2\end{array}\right)\left(\begin{array}{cc}
    \overline{\eta}_2 & \eta_1^{-1} \\ \overline{\eta}_1 & 0
  \end{array}\right) = \left(\begin{array}{cc}
  \overline{\xi}_1^{-1}\overline{\eta}_1 & 0
    \\ \xi_1\overline{\eta}_2+\xi_2\overline{\eta}_1 &
    \xi_1{\eta}_1^{-1}\end{array}\right).$$ This verifies ``$\La$'' and
completes the proof. \ende Finally, as above we can now express this in
terms of ideal classes. The ideals $\aa_{\xi}$ and $\aa_{\eta}$ of
Proposition \ref{p16.1} are equivalent, according to the usual definition.
For general division algebras, one requires a weaker notion of equivalence
than isomorphism as $\Delta$-modules, the notion of {\it stable
  equivalence}. However it is standard that for central simple division
algebras over number fields, excluding definite quaternion algebras, the
stronger notion coincides with the weaker one (see \cite{R}, 35.13).  This
is: $\Delta$-modules $M$ and $N$ are stably equivalent if
$$M+\Delta^r \cong N+\Delta^r,$$ for some $r\geq 0$, ``$\cong$'' being
isomorphism of $\Delta$-modules.

At any rate, one defines the {\it class ray group} (depending on $D$) as
follows. Let $S\inn \Sigma_{\infty}$ denote the set of infinite primes
which ramify in $D$, and set
$$\scC\ell_D({\cal O}_K):=\{ideals\}/\left\{\alpha{\cal O}_K,\ \alpha\in
  K^*, \alpha_{\nu}>0\hbox{ for all } \nu\in S\right\}.$$ Furthermore, let
$\scC\ell(\Delta)$ denote the set of equivalence classes of left ideals in
$\Delta$. The two are related by
\begin{eqnarray*} {\scC}\ell(\Delta) & \stackrel{\sim}{\lra} &
  {\scC}\ell_D({\cal O}_K)
  \\ {[M]} & \mapsto & {[N_{D|K}(M)]}.
\end{eqnarray*}
For details on these matters, see \cite{R}.

Now notice that for the division algebras of (\ref{e1.0}), for $d\geq2$, we
have $\scC\ell_D({\cal O}_K)=\scC\ell({\cal O}_K)$ (where again $K=k$ if
$d=2$).  Hence
\begin{proposition}\label{p16a.1} For the division algebras $D$ we are
  considering, if $d\geq2$ and $\Delta\inn D$ is a maximal order, then
  $$\scC\ell(\Delta)\stackrel{\sim}{\lra}\scC\ell({\cal O}_K).$$
\end{proposition}
We can now proceed to carry out the program sketched above for the division
algebras $D$.
\begin{theorem}\label{t16a.1} Let $D$ be a central simple division algebra
  over $K$ as in (\ref{e1.0}), $D^2$ the right vector space with the
  hyperbolic form (\ref{e4.1}). Let $\xi=(\xi_1,\xi_2)$ and
  $\ge=(\eta_1,\eta_2)$ be two isotropic vectors in $\scI_{\Delta}$. Then
\begin{itemize}\item[(i)] If $d=2$, then $\xi$ and $\eta$ are equivalent
  under $\Gamma_{\Delta}$ $\iff$ the ideals $N_{D|k}(\aa_{\xi})$ and
  $N_{D|k}(\aa_{\eta})$ are equivalent in ${\cal O}_k$.
\item[(ii)] If $d\geq 3$, then $\xi$ and $\eta$ are equivalent under
  $\Gamma_{\Delta}$ $\iff$ the ideals $N_{D|K}(\aa_{\xi})$ and
  $N_{D|K}(\aa_{\eta})$ are equivalent in ${\cal O}_K$.
\end{itemize}
\end{theorem}
{\bf Proof:} By Proposition \ref{p16.1}, $\xi$ and $\eta$ are equivalent
under $\Gamma_{\Delta}$ $\iff$ the ideals $\aa_{\xi}$ and $\aa_{\eta}$ are
equivalent. By Proposition \ref{p16a.1} $\aa_{\xi}\cong \aa_{\eta}$ $\iff
N_{D|K}(\aa_{\xi})\cong N_{D|K}(\aa_{\eta})$ (where $K=k$ for $d=2$). This
is the statement of the theorem. \ende

As a corollary of this
\begin{corollary}\label{c16.1} The number of cusps of $\Gamma_{\Delta}$ is the
  class number of $K$ ($d\geq 3$) or the class number of $k$ ($d=1,2$).
\end{corollary}
{\bf Proof:} It follows from Theorem \ref{t16a.1} that the number of cusps
is {\it at most} the class number in question; it remains to verify that it
is at least the class number. If $\aa_{\xi}$ denotes the ideal in $\Delta$,
and $N_{D|K}(\aa_{\xi})$ denotes the corresponding ideal in ${\cal O}_K$,
then $\xi_1,\ \xi_2$ form a basis of $\aa_{\xi}$, while $N_{D|K}(\xi_1),\
N_{D|K}(\xi_2)$ form an ${\cal O}_K$-basis of the norm ideal. So the
question here is, given an ideal $(a_1,a_2)\inn {\cal O}_K$ generated by
two elements, is there an ideal $(b_1,b_2)$ in the same ideal class of
$(a_1,a_2)$, such that $b_i=N_{D|K}(a_i')$ for some $a_i'\in \Delta$? But
this is just the statement of Proposition \ref{p16a.1}: given any ideal
$\aa$, there is an ideal $\bb\inn[\aa]$, and an ideal $\bb'\inn\Delta$ such
that $\bb=N_{D|K}(\bb')$. We further require that $\bb'$ defines an
isotropic vector in $\Delta^2$, in the above sense. There should be
$b_i\in\bb'$ which generate $\bb'$, such that the vector $(b_1,b_2)\in
\Delta^2$ is isotropic with respect to $h$. This in turn is a question
about what happens to the relation
$\xi_1\overline{\xi}_2+\xi_2\overline{\xi}_1=0$ under the norm map. This
relation turns into
$$N(\xi_1)\overline{N(\xi_2)}+N(\xi_2)\overline{N(\xi_1)}=0,$$ or setting
$b_1=N(\xi_1),\ b_2=N(\xi_2),\ b_1\overline{b}_2+b_2\overline{b}_1=0$. This
just says that $(b_1,b_2)$ is isotropic in the hyperbolic plane $K^2$,
hence $(b_1,b_2)$ is isotropic in the hyperbolic plane $D^2$. In sum, for
any isotropic vector $v\in K^2$, the vector is also an isotropic vector
$v\in D^2$. Any isotropic vector of $D^2$ yields by the norm map an
isotropic vector of $K^2$. It follows from Proposition \ref{p16a.1} that
this gives an isomorphism on equivalence classes, or in other words, {\it
  the isomorphism of Proposition \ref{p16a.1} can be represented
  geometrically by isotropic vectors}. \ende
\begin{remark} We have not been very precise about the groups, but the number
  of cusps of $\Gamma_{\Delta}$ and $S\Gamma_{\Delta}$ are the same, as is
  quite clear.
\end{remark}
Now assume $d\geq 2$ and recall the subfield $L\inn D$ and subgroups
$G_L\inn G_D$ of Proposition \ref{p7.1}. By Proposition \ref{p10.1} these
give rise to subdomains of the domain $\cD_D$. Consider, in $G_L$, the
discrete group $\Gamma_{{\cal O}_L}\inn G_L(K)$.
\begin{lemma}\label{l17.1} We have for any maximal order $\Delta\inn D$,
  $$\Gamma_{\Delta}\cap G_L=\Gamma_{{\cal O}_L}.$$
\end{lemma}
{\bf Proof:} This follows from the definitions and the fact that
$\Delta\cap L={\cal O}_L$. \ende Let $M_L$ denote the arithmetic quotient
\begin{equation}\label{e17.1} M_L=\Gamma_{{\cal O}_L}\bs \cD_L.
\end{equation}
It follows from Proposition \ref{p10.1} that we have a commutative diagram
\begin{equation}\label{e17.2} \begin{array}{rcl} \cD_L & \hra & \cD_D \\
    \downarrow & & \downarrow \\ M_L & \hra & X_{\Gamma_{\Delta}}.
\end{array}\end{equation}
Each $M_L$ has its own Baily-Borel compactification $M_L^*$, which is also
affected by adding isolated points to $M_L$. Moreover, from Proposition
\ref{p10.2}, we see that the embedding $M_L\hra X_{\Gamma_{\Delta}}$ can be
extended to the cusp denoted $F_L\in \cD_L^*$ respectively $F\in \cD_D^*$
there. By Theorem 3 of \cite{S2}, we actually get embeddings of the
Baily-Borel embeddings of $M_L$ and $X_{\Gamma_{\Delta}}$, respectively.
\begin{theorem}\label{t17.1} Let $M_L^*\inn {\Bbb P}^N,\
  X_{\Gamma_{\Delta}}^*\inn
  {\Bbb P}^{N'}$ be Baily-Borel embeddings. Then there is a linear
  injection ${\Bbb P}^N\hra {\Bbb P}^{N'}$ making the diagram
  $$\begin{array}{rcl} M_L^* & \hra & {\Bbb P}^N \\ \cap & & \cap \\
    X_{\Gamma_{\Delta}}^* & \hra & {\Bbb P}^{N'}\end{array}$$ commute and
  making $M_L^*\in X_{\Gamma_{\Delta}}^*$ an algebraic subvariety.
\end{theorem}
{\bf Proof:} We have an injective holomorphic embedding $\cD_L\hra \cD_D$
which comes from a ${\Bbb Q}$-morphism $\rho:(G_L)_{{\Bbb C}}\hra
(G_D)_{{\Bbb C}}$ (for this one takes the restriction of scalars from $k$
to ${\Bbb Q}$ yielding an injection $Res_{k|{\Bbb Q}}G_L\hra Res_{k|{\Bbb
    Q}}G_D$, then lifts this to ${\Bbb C}$) such that $\rho(\Gamma_{{\cal
    O}_L})\inn \Gamma_{\Delta}$. Hence we map apply \cite{S2}, Thm.~3, and
the theorem follows. \ende
\begin{corollary}\label{c17.1} If $d=2$, then there are modular curves
  $M_L^*$ on the algebraic threefold $X_{\Gamma_{\Delta}}^*$ such that the
  cusps of $M_L^*$ are cusps of $X_{\Gamma_{\Delta}}^*$. If $d\geq 3$, then
  we have Hilbert modular varieties of dimension $d$, $M_L^*\inn
  X_{\Gamma_{\Delta}}^*$ in the $d^2$-dimensional algebraic variety
  $X_{\Gamma_{\Delta}}^*$, such that the cusps of $M_L^*$ are cusps of
  $X_{\Gamma_{\Delta}}^*$.
\end{corollary}
The previous Theorem \ref{t17.1} applies to the cusp of
$X_{\Gamma_{\Delta}}^*$ which represents the equivalence class of the
isotropic vector $(0,1)$. We now consider the others. Given
$(\xi_1,\xi_2)$, an isotropic vector representing a class of cusps, let
$M_{\xi}$ be the matrix in $G_D$ of Lemma \ref{l15.2} which maps it to
$(0,1)$. Then
$$\Gamma_{\Delta,\xi}:=M_{\xi}\cdot \Gamma_{\Delta}\cdot M_{\xi}^{-1}\inn
U(D^2,h)$$ is a discrete subgroup of $G_D$, and the cusp $(\xi_1,\xi_2)$ of
$\Gamma_{\Delta}$ is equivalent to the cusp $(0,1)$ of
$\Gamma_{\Delta,\xi}$. Letting, as above, $\Gamma_{{\cal O}_L}\inn
\Gamma_{\Delta}$ denote the discrete subgroup defining the modular
subvariety $M_L$ above, we have
\begin{equation}\label{e17a.0} \Gamma_{{\cal O}_L,\xi}:=M_{\xi}\cdot
  \Gamma_{{\cal O}_L}\cdot M_{\xi}^{-1}\inn \Gamma_{\Delta,\xi},
\end{equation}
and without difficulty this gives a subdomain
\begin{equation}\label{e17a.1} \cD_{L,\xi}\inn \cD_D,
\end{equation}
such that, if $F_{\xi}$ denotes the boundary component corresponding to
$\xi$, then $F_{\xi}\in \cD_{L,\xi}^*$. More precisely,
$G_{L,\xi}:=M_{\xi}\cdot G_L\cdot M_{\xi}^{-1}$ is the $k$-group whose
${\Bbb R}$-points $G_{L,\xi}({\Bbb R})$ define $\cD_{L,\xi}$. Then the
parabolic subgroup of $\Gamma_{\Delta,\xi}$ is $\Gamma_{\Delta,\xi}\cap P$
($P$ as in (\ref{e6.1}) the normaliser of $(0,1)$), and similarly for
$\Gamma_{{\cal O}_L,\xi}$. The corresponding modular subvariety
\begin{equation}\label{e17a.2} M_{L,\xi}=\Gamma_{{\cal O}_L,\xi}\bs
  \cD_{L,\xi},
\end{equation}
is, as is easily checked, a Hilbert modular variety for the group
$SL_2({\cal O}_{\ell},\bb^2)$, where for any ideal $\cc$ one defines
\begin{equation}\label{e17a.3}
  SL_2({\cal O}_{\ell},\cc)=\left\{\left(\begin{array}{cc} a & b \\ c &
        d\end{array}\right) \Big| ad-cb=1,\ a,d\in {\cal O}_{\ell}, b\in
    \cc^{-1},\ c\in \cc\right\},\end{equation} and where $\bb$ is the
intersection of the ideal $\aa_{\xi}$ with ${\cal O}_{\ell}$. The same
arguments as above then yield
\begin{theorem}\label{t17a.1} Given any cusp $p\in X_{\Gamma_{\Delta}}^*\bs
  X_{\Gamma_{\Delta}}$, there is a modular subvariety $M_{L,p}\inn
  X_{\Gamma_{\Delta}}$ such that $p\in M_{L,p}^*$.
\end{theorem}
{\bf Proof:} If the cusp $p$ is represented by the isotropic vector
$\xi=(\xi_1,\xi_2)$, then the modular subvariety is $M_{L,\xi}$ as in
(\ref{e17a.2}). Just as in the proof of Theorem \ref{t17.1}, we get an
embedding of the Baily-Borel embeddings, and as mentioned above,
$F_{\xi}\in \cD_{L,\xi}^*$, and $p$ is the image of $F_{\xi}$ under the
natural projection $\pi:\cD_{L,\xi}^*\lra M_{L,\xi}^*=M_{L,p}^*$. \ende
\begin{corollary}\label{c17a.1} If $d=2$, there are modular curves
  $M_{L,p}^*\inn X_{\Gamma_{\Delta}}^*$ passing through each cusp of
  $X_{\Gamma_{\Delta}}^*$. If $d\geq3$, there are Hilbert modular varieties
  (to groups as in (\ref{e17a.3})) of dimension $d$ passing through any
  cusp of $X_{\Gamma_{\Delta}}^*$. \end{corollary}

\section{An example}
The theory of quaternion algebras is quite well established, and the
corresponding arithmetic quotients have already been studied in \cite{Ara}
and \cite{H}. However the $d\geq3$ case seems not to have drawn much
attention as of yet, so we will give an example to illustrate the theory.
Curiously enough, I ran across this example in the construction of
Mumford's fake projective plane. Recall that this is an algebraic surface
$S$ of general type with $c_1^2=3c_2,\ c_1^2=9, \ c_2=3$, just as for the
projective plane. It then follows from Yau's theorem that $S$ is the
quotient of the two-ball ${\Bbb B}_2$ by a discrete subgroup,
\begin{equation}\label{e18.1} S=\Gamma\bs {\Bbb B}_2,
\end{equation}
where $\Gamma$ is cocompact and fixed point free. Mumford's construction
involves lifting a quotient surface from a $2$-adic field to the complex
numbers, and while the $2$-adic group is clearly not arithmetic (as a
subgroup of $SL(3,{\Bbb Q}_2)$, its elements are unbounded in the $2$-adic
valuation), it is not clear from the construction whether $\Gamma$ is
arithmetic. Now the {\it arithmetic} cocompact groups are known: these
derive either from anisotropic groups (over ${\Bbb Q}$) which are of the
form $U(1,D)$ or $SU(1,D)$, where $D$ is a central simple division algebra
over an imaginary quadratic extension $K$ of a totally real field $k$, such
that if the corresponding hermitian symmetric domain is irreducible, then
$k={\Bbb Q}$, and such that $D$ has a $K|{\Bbb Q}$-involution, or from
unitary groups over field extensions $k|{\Bbb Q}$ of degree $d\geq2$, such
that, for all but one infinite prime, the real groups $G_{\nu}$ are
compact. The strange thing is that in Mumford's construction such a $D$
comes up (implicitly) naturally, and it is this $D$ we will introduce. It
is a fascinating question whether the $\Gamma$ of (\ref{e18.1}) occurs as a
discrete subgroup of $U(1,D)$, a question which will be left unanswered
here.

This is a cyclic algebra, central simple over $K$, with splitting field
$L$. The fields involved are
\begin{equation}\label{e18.2} L={\Bbb Q}(\zeta),\ \ \zeta=\exp({2\pi i \over
    7});\quad K={\Bbb Q}(\sqrt{-7}).
\end{equation}
\begin{lemma}\label{l18.1} $L$ is a cyclic extension of $K$, of degree
  three. The Galois group is generated by the transformation
  $$\sigma(\zeta)=\zeta^2,\ \sigma(\zeta^2)=\zeta^4,\
  \sigma(\zeta^4)=\zeta,$$ where $1, \zeta, \zeta^2$ generate $L$ over $K$.
\end{lemma}
{\bf Proof:} Let $\gamma={-1+\sqrt{-7}\over 2}\in K$, so ${\cal O}_K={\Bbb
  Z}\oplus \gamma{\Bbb Z}$. The ring of integers of $L$ is generated by the
powers of $\zeta$,
$${\cal O}_L={\Bbb Z}\oplus\zeta{\Bbb Z}\oplus \cdots \oplus \zeta^5{\Bbb
  Z},$$ (remember that $\zeta^6=-1-\zeta-\cdots-\zeta^5$), with the
inclusion ${\Bbb Z}\inn {\cal O}_K, {\Bbb Z}\inn {\cal O}_L$ on the first
factors. The key identity is the following
\begin{equation}\label{e18.3} \gamma=\zeta+\zeta^2+\zeta^4,\quad\quad
  \overline{\gamma}=\zeta^3+\zeta^5+\zeta^6.
\end{equation}
{}From this relation we see that ${\cal O}_L$ is generated over ${\cal O}_K$
by $\zeta$ and $\zeta^2$,
$${\cal O}_L={\cal O}_K\oplus\zeta{\cal O}_K\oplus\zeta^2{\cal O}_K.$$ The
conjugation on $K$ extends to (complex) conjugation on $L$, which we will
denote by $x\mapsto \overline{x}$. Its action on ${\cal O}_K$ is clear,
$\gamma\mapsto \overline{\gamma}$, and on ${\cal O}_L$ it affects
$$\overline{\zeta}=\zeta^6,\quad \overline{\zeta}^2=\zeta^5,\quad
\overline{\zeta}^3=\zeta^4.$$ It is clear that $L|K$ is cyclic, with a
generator of the Galois group
$$1\mapsto 1, \quad \zeta\mapsto \zeta^2,\quad \zeta^2\mapsto
\gamma-\zeta-\zeta^2,$$ and by (\ref{e18.3}), this is the statement of the
Lemma. \ende Now consider the fixed field in $L$ under conjugation,
$\ell\inn L$.
\begin{lemma}\label{l19.1} $\ell$ is a degree three cyclic extension of
  ${\Bbb Q}$, with a generator of the Galois group being
  $$\sigma:\quad \zeta+\zeta^6\mapsto \zeta^2+\zeta^5 \mapsto
  \zeta^3+\zeta^4\mapsto \zeta+\zeta^6.$$ \end{lemma} {\bf Proof:} It is
clear that $\eta_1=\zeta+\zeta^6, \eta_2=\zeta^2+\zeta^5,\
\eta_3=\zeta^3+\zeta^4$ are in $\ell$ (they are
$\eta_1=Tr_{L|\ell}(\zeta),\ \eta_2=Tr_{L|\ell}(\zeta^2),\
\eta_3=Tr_{L|\ell}(\zeta^3)$).  Moreover they generate the integral closure
of ${\Bbb Z}$ in $\ell$, hence $${\cal O}_{\ell}={\Bbb Z}\oplus \eta_1{\Bbb
  Z}\oplus \eta_2{\Bbb Z},$$ and there is a relation
$\eta_3=-\eta_1-\eta_2-1$.  Now observe that $\sigma$ of Lemma \ref{l18.1}
acts as follows
$$\sigma(\eta_1)=\sigma(\zeta+\zeta^6)=\sigma(\zeta)+\sigma(\zeta^6)=
  \zeta^2+\zeta^5=\eta_2,\quad
\sigma(\eta_2)=\eta_3, \quad \sigma(\eta_3)=\eta_1,$$ and this is the map
specified in the lemma. \ende Next consider the cyclic algebra
$D=(L/K,\sigma,\gamma)$ for the element $\gamma$ above. Then $D$ will be
split $\iff$ $\gamma$ is the $L/K$ norm of some element. However, it seems
difficult to verify this condition explicitly, so we show directly that $D$
is a division algebra.
\begin{proposition}\label{p19.1} $D=(L/K,\sigma,\gamma)$, with
  $\gamma={-1+\sqrt{-7}\over 2}$, is a division algebra, central simple
  over $K$.
\end{proposition}
{\bf Proof:} Note first that since the degree of $D$ over $K$ is three, a
prime, $D$ is a division algebra $\iff$ $D$ is not split. To show that $D$
is not split, it suffices to find a prime $\pp\in {\cal O}_K$ for which the
local algebra $D_{\pp}$ is not split. This, it turns out, is easy to find.
Quite generally we know that $D_{\pp}$ is split for almost all $\pp$, and
non-split at the divisors of the discriminant. Note that
$$N_{K|{\Bbb Q}}(\gamma)=\gamma\overline{\gamma}=2,$$ while $K$ ramifies
over ${\Bbb Q}$ only at the prime $\sqrt{-7}$. These are the two primes
where $D_{\pp}$ might ramify.  To see that for $\pp=(\gamma)$ $D_{\pp}$
actually does ramify, we determine the Hasse invariant $inv_{\pp}(D)\in
{\Bbb Q}/{\Bbb Z}$ of $D_{\pp}$.

To do this we must first understand the action of the Frobenius acting as
generator for the maximal unramified extension $L_{\pp}/K_{\pp}$. But
Frobenius in characteristic two is just a squaring map, so is clearly the
$\mod(\pp)$ reduction of the $\sigma$ of Lemma \ref{l18.1}. We denote this
by $\sigma_{\pp}$. Finally, since we are at the prime $\pp=(\gamma)$, we
may take the image of $\gamma$ as local uniformising element $\pi_{\pp}$.
Then
$$D_{\pp}=(L_{\pp}/K_{\pp},\sigma_{\pp},\pi_{\pp}),$$ which is the cyclic
algebra with $inv_{\pp}(D)={1\over 3}$. From this it follows that $D_{\pp}$
is not split. That $D$ is central simple over $K$ is clear from
construction. \ende Finally we require a $K|{\Bbb Q}$-involution on $D$. It
is necessary and sufficient for the existence of such an involution $J$
that there exists an element $g\in \ell$, such that
\begin{equation}\label{e20.1} N_{K|{\Bbb Q}}(\gamma)=
  \gamma\overline{\gamma}=N_{\ell|{\Bbb Q}}(g) =g
  g^{\sigma}g^{\sigma^2}.
\end{equation}
Given such an element $g$, the involution is given by (see (\ref{e4A.1}))
$$e^J=ge^{-1},\ \ (e^2)^J=g g^{\sigma}(e^2)^{-1},\quad (\sum_0^2e^iz_i)^J:=
\sum_0^2\overline{z}_i(e^i)^J.$$ An alternative to finding an explicit $g$
is to use a theorem of Landherr (\cite{Sch}, Thm.~10.2.4): $D$ admits a
$K|{\Bbb Q}$-involution $\iff$
\begin{equation}\label{e20.2} \begin{minipage}{15cm} \begin{itemize}\item
      $inv_{\pp}(D)=0,\ \quad \forall_{\pp=\overline{\pp}}$, and
    \item $inv_{\pp}(D)+inv_{\overline{\pp}}(D)=0, \quad \forall_{\pp\neq
        \overline{\pp}}$.
\end{itemize}\end{minipage}\end{equation}
This condition is satisfied for all $\pp$ for which $D_{\pp}$ splits, so
must be verified only for those primes which divide the discriminant. For
us, this means at most $\pm\sqrt{-7}, \gamma,\ \overline{\gamma}$, which we
will denote by $\pm \qq, \pp, \overline{\pp}$, respectively. We showed
above that $inv_{\pp}(D)={1\over 3}$. The same argument shows that
$inv_{\overline{\pp}}(D)=-{1\over 3}$. The negative sign occurs because at
the prime $\overline{\pp}$, $\overline{\gamma}$ is a local uniformising
element $\pi_{\overline{\pp}}$, so localising $\gamma$ at $\overline{\pp}$
gives $\pi_{\overline{\pp}}^{-1}$. This verifies (\ref{e20.2}) for the
primes $\pp\neq\overline{\pp}$ lying over 2.  Consider the prime $\pm \qq$.
Here it is easy to see that the image of $\gamma$ is actually {\it
  invertible}, and hence, that $D_{\qq}$ splits, so $inv_{\qq}(D)=0$, as
was to be shown.  \ende This completes the proof of the following.
\begin{theorem}\label{t20.1} The cyclic algebra $D=(L/K,\sigma,\gamma)$
  constructed above is a central simple division algebra over $K$ with a
  $K|{\Bbb Q}$-involution. It ramifies at the two primes $\gamma$ and
  $\overline{\gamma}$, and is split at all others.
\end{theorem}
This algebra gives rise to an anisotropic group $GL(1,D)\cong D^*$, and
$D^*({\Bbb R})$ is a twisted real form of $GL_3$; since it cannot be
compact, it must be $U(2,1)$. A maximal order $\Delta\inn D$ gives rise to
an arithmetic subgroup of $D^*({\Bbb R})$, and the quotient is then a
compact ball quotient.  As remarked above, this may be related with
Mumford's fake projective plane.

We can consider the hyperbolic space $(D^2,h)$, and the corresponding
${\Bbb Q}$-groups (here $k={\Bbb Q}$) $G_D$ and $SG_D$ as in (\ref{e4.2})
and (\ref{e4.4}), respectively, as well as the arithmetic subgroups
$\Gamma_{\Delta}$ and $S\Gamma_{\Delta}$. By Corollary \ref{c16.1}, we see
that $\Gamma_{\Delta}$ has $h(K)$ cusps, where $h(K)$ is the class number
of $K$; this is known to be 1. So the arithmetic quotient has only 1 cusp.
We also have the subgroups $\Gamma_{{\cal O}_L}$ as in Lemma \ref{l17.1},
as well as the modular subvarieties $M_L^*$ of $X_{\Gamma_{\Delta}}^*$. We
note that these are Hilbert modular threefolds coming from the cyclic cubic
extension $\ell/{\Bbb Q}$. Such threefolds have been considered in
\cite{T}.

\section{Moduli interpretation}
The moduli interpretation of the arithmetic quotients
$X_{\Gamma_{\Delta}}$, as well as of the modular subvarieties $M_L$, by
which we mean the description of these spaces as moduli spaces, is a
straightforward application of Shimura's theory. Fix a hyperbolic plane
$(D^2,h)$, and consider the endomorphism algebra $M_2(D)$ of $D^2$, endowed
with the involution $x\mapsto {^t\overline{x}}$, $x\in M_2(D)$, where the
bar denotes the involution in $D$. Now view $D^2$ as a ${\Bbb Q}$-vector
space, of dimension $4f$, $8f$ and $4d^2f$, in the cases $d=1$, $d=2$ and
$d\geq3$, respectively. The data determining one of Shimura's moduli spaces
(with no level structure) is $(D,\Phi,*),\ (T,{\cal M})$, where $D$ is a
central simple division algebra over $K$, $\Phi$ is a representation of $D$
in $\Gg\ll(n,{\Bbb C})$, $*$ is an involution on $D$, $T$ is a $*$-skew
hermitian form (matrix) on a right $D$ vector space $V$ of dimension $m$,
with $\Gg\ll(n,{\Bbb C})\cong End(V,V)_{{\Bbb R}}$, and finally, ${\cal
  M}\inn D^m$ is a ${\Bbb Z}$-lattice. Then for suitable $x_i\in V$,
\begin{equation}\label{e21.1} \Lambda=\left\{\sum_{i=1}^m\Phi(a_i)x_i
  \Big| a_i\in
    {\cal M}\right\}\inn {\Bbb C}^n
\end{equation}
is a lattice and ${\Bbb C}^n/\Lambda$ is abelian variety with
multiplication by $D$.  The data $(D,\Phi,*)$ will be given in our cases as
follows. $D$ is our central simple division algebra over $K$ with a
$K|k$-involution, and the representation $\Phi:D\lra M_N({\Bbb C})$ is
obtained by base change from the natural operation of $D$ on $D^2$ by right
multiplication. Explicitly,
\begin{equation}\label{e21.3} \Phi:D\lra End_D(D^2,D^2)
  \otimes_{{\Bbb Q}}{\Bbb R}
  \cong M_2(D)\otimes_{{\Bbb Q}}{\Bbb R}\cong M_2(D\otimes_{{\Bbb Q}}{\Bbb
    R})\cong M_2({\Bbb R}^N)\cong M_N({\Bbb C}),
\end{equation}
where $N=\dim_{{\Bbb Q}}D=2f,\ 4f,\ 2d^2f$ in the cases $d=1,\ d=2$ and
$d\geq3$, respectively. The involution $*$ on $D$ will be our involution,
which we still denote by $x\mapsto \overline{x}$. Then a $^{-}$-skew
hermitian matrix $T\in M_2(D)$ will be one such that $T=-T^*$, where
$(t_{ij})^*=(\overline{t}_{ji})$, the canonical involution on $M_2(D)$
induced by the involution on $D$. Note that for any $c\in D^*$ such that
$c=-\overline{c}$, the matrix $T=cH$ ($H$ our hyperbolic matrix ${0\
  1\choose 1\ 0}$) has this property. To be more specific, then, we set
\begin{equation}\label{e21.4}
\begin{minipage}{15cm}\begin{itemize}\item[1)] $d=1$:
    $T=\sqrt{-\eta}H=\left(\begin{array}{cc} 0 & \sqrt{-\eta} \\
        \sqrt{-\eta} & 0 \end{array}\right)$.
  \item[2)] $d=2$: $T_a=\sqrt{a}H=\left(\begin{array}{cc} 0 & \sqrt{a} \\
        \sqrt{a} & 0 \end{array}\right)$ or $T_b=eH={0\ e\choose e\ 0}$,
    where $e={0\ 1\choose b\ 0}$.
  \item[3)] $d\geq3$: $T=\sqrt{-\eta}H=\left(\begin{array}{cc} 0 &
        \sqrt{-\eta} \\ \sqrt{-\eta} & 0 \end{array}\right)$.
\end{itemize}\end{minipage}\end{equation}
Since two such forms $T$ are equivalent when they are scalar multiples of
one another, assuming $T$ of the form in (\ref{e21.4}) is no real
restriction. Finally the lattice ${\cal M}$ will be $\Delta^2\inn D^2$.
Then $D^2\otimes_{{\Bbb Q}}{\Bbb R}\cong {\Bbb C}^N$, $N$ as above, and for
``suitable'' vectors $x_1, x_2\in D^2$, the lattice
\begin{equation}\label{e21.2} \Lambda_x=\{\Phi(a_1)x_1+\Phi(a_2)x_2 \big|
  (a_1,a_2)\in \Delta^2\}
\end{equation}
gives rise to an abelian variety $A_x={\Bbb C}^N/\Lambda_x$. Shimura has
determined exactly what ``suitable'' means; the conditions determine
certain unbounded realisations of hermitian symmetric spaces, in our cases
just the domains $\cD_D$. The Riemann form on $A_x$ is given by the
alternating form $E(x,y)$ on ${\Bbb C}^N$ defined by:
\begin{equation}\label{ee21.6} E(\sum_1^2\Phi(\alpha_i)x_i,
  \sum_1^2\Phi(\beta_j)x_j)=Tr_{D|{\Bbb Q}}(\sum_{i,j=1}^2
  \alpha_it_{ij}\overline{\beta}_j),
\end{equation}
for $\alpha_i,\beta_j\in D_{{\Bbb R}}$, and $(t_{ij})=T$ is the matrix
above. In particular, in all cases dealt with here the abelian varieties
are {\it principally polarised}.

Shimura shows that for each $z\in \cD_D$, vectors $x_1,x_2\in D^2$ are
uniquely determined, hence by (\ref{e21.2}) a lattice, denoted
$\Lambda(z,T,{\cal M})$. The data determine an arithmetic group, which for
our cases is just $\Gamma_{\Delta}$ defined above (not $S\Gamma_{\Delta}$),
cf.~\cite{Sh2} (38).  The basic result, applied to our concrete situation,
is
\begin{theorem}[\cite{Sh2}, Thm.~2]\label{t21.1} The arithmetic quotient
  $X_{\Gamma_{\Delta}}$ is the moduli space of isomorphism classes of
  abelian varieties determined by the data:
  $$(D,\Phi,^-),\ \ (T,\Delta^2),$$ where $\Phi$ is given in (\ref{e21.3}),
  $T$ in (\ref{e21.4}).
\end{theorem}
The corresponding classes of abelian varieties can be described as follows:
\begin{itemize}\item[1)] $d=1$. Here we have {\it two} families, relating
  from the isomorphism of Proposition \ref{p5.1}. The first, for $D=k,\
  D^2=k^2$ yields $D_{{\Bbb R}}^2\cong {\Bbb C}^f$, and we have abelian
  varieties of dimension $f$ with real multiplication by $k$. Secondly, for
  $D=K,\ D^2=K^2$, we have abelian varieties of dimension $2f$ with complex
  mulitplication by $K$, with signature $(1,1)$, that is, for each
  eigenvalue $\chi$ of the differential of the action, also
  $\overline{\chi}$ occurs. If $K=k(\sqrt{-\eta})$, then setting $K'={\Bbb
    Q}(\sqrt{-\eta})$ we have $k\otimes_{{\Bbb Q}}K'\cong K$, hence
  $k^2\otimes_{{\Bbb Q}}K'\cong K^2$ and $k^2_{{\Bbb R}}\otimes K'_{{\Bbb
      R}}\cong K_{{\Bbb R}}^2$, giving the relation between the ${\Bbb
    Q}$-vector spaces and their real points. Moreover, ${\cal
    O}_k\otimes_{{\Bbb Z}}{\cal O}_{K'}\cong {\cal O}_K$, and if
\begin{equation}\label{e21.6} \Lambda_{x,k}=\{\sum_1^2\Phi(a_i)x_i \Big|
  (a_1,a_2)\in {\cal O}_k^2\}
\end{equation}
is a lattice giving an abelian variety with mulitiplication by $k$,
$$A_x:={\Bbb C}^f/\Lambda_{x,k},$$ then $\Lambda_{x,k}\otimes_{{\Bbb
    Z}}{\cal O}_{K'}=\Lambda_{x,K}$ is a lattice in ${\Bbb C}^{2f}$, and
determines an abelian variety
\begin{equation}\label{e21.7} A_x':={\Bbb C}^{2f}/\Lambda_{x,K}.
\end{equation}
This abelian variety determines a point $x'$ in its moduli space, and the
mapping $x\mapsto x'$ gives the isomorphism
\begin{equation}\label{e21.5} \Gamma_{{\cal O}_k}\bs \HH^f
  \stackrel{\sim}{\lra}
  \Gamma_{{\cal O}_K}\bs \HH^{f}.
\end{equation}
\begin{remark} It turns out that this case is one of the exceptions of
  Theorem 5 in \cite{Sh2}, denoted case d) there. The actual endomorphism
  ring of the generic member of the family is larger than $K$:
\begin{theorem}[\cite{Sh2}, Prop.~18] The endomorphism ring $E$ of the
  generic element of the family (\ref{e21.5}) is a totally indefinite
  quaternion algebra over $k$, having $K$ as a quadratic subfield.
\end{theorem}
In our situation, the totally indefinite quaternion algebra $E$ over $k$ is
constructed as the cyclic algebra $E=(K/k,\sigma,\lambda)$, where
$\lambda=-u^{-1}v$, if the matrix $T$ of (\ref{e21.4}) is diagonalized
$T={u\ 0\choose 0\ v}$.  So in our case we have $\lambda=1$ and hence the
algebra $E$ is split; the corresponding abelian variety is isogenous to a
product of two copies of a simple abelian variety $B$ with real
multiplication by $k$, as has been described already above. The conclusion
follows from our choice of $T$, i.e., of hyperbolic form. It would seem one
gets more interesting quaternion algebras by choosing different hermitian
forms (which, by the way, will also lead to other polarisations).
\end{remark}
\item[2)] $d=2$. $D=(\ell/k,\sigma,b)=(a,b)$ is a totally indefinite
  quaternion algebra, cental simple over $k$, with canonical involution. We
  have $D_{\nu}\cong M_2({\Bbb R})$, and $D\otimes_{{\Bbb Q}}{\Bbb R}\cong
  {\Bbb R}^4$, while $M_2(D)\otimes_{{\Bbb Q}}{\Bbb R}\cong M_2({\Bbb
    R}^4)\cong M_4({\Bbb C})$. Let $\Delta\inn D$ be a maximal order,
  $\Gamma_{\Delta} \inn G_D$ the corresponding arithmetic group.  Two
  vectors $x_1,x_2\in D^2$ arising from a point in the domain ${\Bbb S}_2$
  (Siegel space of degree 2) determine a lattice $\Lambda_x$ as in
  (\ref{e21.6}), with $(a_1,a_2)\in \Delta^2$, and $A_x={\Bbb
    C}^{4f}/\Lambda_x$ is the corresponding abelian variety.
\item[3)] $d\geq3$. In this case $D$ is the cyclic algebra of degree $d$
  over $K$, and the abelian varieties are of dimension $2d^2f$.
\end{itemize}

Now we come to the most interesting point of the whole story -- the moduli
interpretation of the modular subvarieties $M_L\inn X_{\Gamma_{\Delta}}$ of
(\ref{e17.2}). The moduli interpretation of the $M_L$ is as stated in
Theorem \ref{t21.1}, $d=1$ case. Disregarding the case $d=1$, the
subvarieties have the following interpretations.
\begin{itemize}\item[2)] $d=2$:
  As $M_L$ arises from the group $U(L^2,h)$, where $L=k(\sqrt{-ab})$ in our
  notations above, this is the moduli space of abelian varieties of
  dimension $2f$ with complex multiplication by $L$. On the other hand, the
  space $X_{\Gamma_{\Delta}}$ parameterizes abelian varieties of dimension
  $4f$ with multiplication by $D$. The relation is given as follows. By
  definition we have $D=\ell\oplus e\ell$, which in terms of matrices, is
  $$D\cong \left\{\left(\begin{array}{cc} a_0+a_1\sqrt{a} & 0 \\ 0 &
        a_0-a_1\sqrt{a}\end{array}\right) \oplus\left(\begin{array}{cc} 0 &
        a_2+a_3\sqrt{a} \\ b(a_2-a_3\sqrt{a}) & 0
      \end{array}\right)\right\}.$$ Now our subfield $L$ is the set of
  matrices of the form
  $$L\cong \left\{\left(\begin{array}{cc} a_0 & a_3\sqrt{a} \\
        -ba_3\sqrt{a} & a_0\end{array}\right)\right\}.$$ If we let
  $c=\diag(\sqrt{a},-\sqrt{a})$ be the element representing $\sqrt{a}$ and
  $e={0\ 1\choose b\ 0}$, then $c(ec)=c(-ce)=-c^2e=-ae$, so we can generate
  $D$ over ${\Bbb Q}$ by $c$ and $(ec)$. Recall that $L\cong k(ec)$, hence
\begin{equation}\label{e22a.1} D\cong L\oplus cL.\end{equation}
Now consider the representation $\Phi$; we have $\Phi(D)=\Phi(L\oplus cL) =
\Phi_{|L}(L)\oplus \Phi_{|cL}(cL)$. The lattice $\Delta^2\inn D^2$ gives
rise to a lattice $\Lambda_x$ as in (\ref{e21.2}), and we would like to
determine when the splitting (\ref{e22a.1}) gives rise to a splitting of
the lattice $\Lambda_x$, hence of the abelian variety $A_x$. Consider the
order
\begin{equation}\label{e22a.2} \Delta':={\cal O}_L\oplus c{\cal O}_L\inn
  \Delta;
\end{equation}
$\Delta'$ is in general not a maximal order, but it is of finite index in
$\Delta$. Note that $\Delta'$ and a point $x$ (consisting of two vectors
$x_1, x_2\in D^2$) determine a lattice
$$\Lambda'_x=\{\sum \Phi(a_i)x_i \big| (a_1,a_2)\in \Delta'\},$$ which is
also of finite index in $\Lambda_x$. Therefore $A'_x={\Bbb C}^4/\Lambda_x'$
and $A_x$ are isogenous.

We now assume $x_i\in L^2$. From (\ref{e22a.2}) we can write
$a_i=a_i^1+ca_i^2$ for $a_i\in \Delta'$, hence
$\Phi(a_i)=\Phi(a_i^1+ca_i^2)= \Phi_{|{\cal O}_L}(a_i^1)+c\Phi_{|{\cal
    O}_L}(a_i^2). $ Then we have
\begin{eqnarray}\label{e22.1} \Lambda_x'=\left\{\sum\Phi(a_i)x_i \Big|
  (a_1,a_2)\in \Delta'^2\right\} & = &
  \left\{\sum\left(\Phi_{|L}(a_i^1)+c\Phi_{|L}(a_i^2)\right) x_i \Big|
    (a_1,a_2)\in \Delta'^2 \right\} \\ \nonumber & = & \left\{
    \Phi_{|L}(a_1^1)x_1+\Phi_{|L}(a_2^1)x_2+c\left(\Phi_{|L}(a_1^2)x_1 +
      \Phi_{|L}(a_2^2)x_2\right) \right\} \\ \nonumber & = &
  \Lambda^1\oplus c \Lambda^2,\nonumber
\end{eqnarray}
and each of $\Lambda^i$ has complex multiplication by $L$. It follows from
this that
$$A_x'\cong A_x'^1\times A_x'^2,\quad x=(x_1,x_2)\in L^2,$$ and each
abelian variety $A_x'^i$, of dimension $2f$, has complex multiplication by
$L$. Since $A_x'\lra A_x$ is an isogeny, we have
\begin{proposition}\label{p22.1} In case $d=2$, the abelian varieties
  parameterised by the modular subvariety $M_L$ are isogenous to products
  of two abelian varieties of dimensions $2f$ with complex multiplication
  by $L$. \end{proposition}
\item[3)] $d\geq3$: Again the situation is somewhat simpler here. As above,
  the subvariety $M_L$ parameterises abelian varieties of dimension $2df$
  with complex multiplication by $L$.  Again we consider the order
  $$\Delta':={\cal O}_L\oplus e{\cal O}_L\oplus \cdots \oplus e^{d-1}{\cal
    O}_L\inn\Delta,$$ where we assume $e$ is as in (\ref{e1.2}) and
  $\gamma\in {\cal O}_k$; then as above we can write the lattice
  $\Lambda'_x$ for $x_1, x_2\in L^2$, in terms of $\Phi_{|L}$. If we write
  $a_i=a_i^1+ea_i^2+\cdots +e^{d-1}a_i^d$ with $a_i^j\in {\cal O}_L$, then
  $\Phi(a_i)=\Phi_{|L}(a_i^1)+e\Phi_{|L}(a_i^2)+\cdots+e^{d-1}
  \Phi_{|L}(a_i^d)$, and consequently
\begin{eqnarray*}\Phi(a_1)x_1+\Phi(a_2)x_2 & = &
  \left(\Phi_{|L}(a_1^1)x_1+\Phi_{|L}(a_2^1)x_2\right) +
  e\left(\Phi_{|L}(a_1^2)x_1+\Phi_{|L}(a_2^2)x_2\right) + \\ & & \ \ +
  \ldots + e^{d-1}\left(\Phi_{|L}(a_1^d)x_1+\Phi_{|L}(a_2^d)x_2\right), \\
  \left\{\sum\Phi(a_i)x_i\Big| (a_1,a_2)\in {\Delta'}^2\right\} & = &
  \Lambda^1\oplus e\Lambda^2\oplus \cdots \oplus e^{d-1}\Lambda^d,
\end{eqnarray*} and each sublattice $\Lambda^i$ has complex multiplication
by $L$. Again, the abelian variety $A_x'$ so determined is isogenous to
$A_x$, and hence
\begin{proposition}\label{p22.2} In case $d\geq3$, the abelian varieties
  parameterised by the modular subvariety $M_L$ are isogenous to the
  product of $d$ abelian varieties of dimension $2df$ with complex
  multiplication by the field $L$.
\end{proposition}
\end{itemize}


\begin{thebibliography}{MMM}

\bibitem[A]{A} A. A. Albert, ``Structure of algebras'', AMS Colloqium Pub.
  XXIV, AMS: Providence 1961.

\bibitem[Ara]{Ara} T. Arakawa, The dimension of the space of cusp forms on
  the Siegel upper half plane of degree two related to a quaternion unitary
  group, J.~Math.~Soc.~Japan {\bf 33} (1981) 125-145.

\bibitem[Ar]{Ar} E. Artin, ``Geometric algebra'', Wiley Interscience: New
  York 1957.

\bibitem[H]{H} K. Hashimoto, The dimension of the spaces of cusp forms on
  Siegel upper half-plane of degree two, Math. Ann. {\bf 266} (1984),
  539-559.

\bibitem[R]{R} I. Reiner, ``Maximal orders'', Acedemic Press: New York
  1969.

\bibitem[S1]{S} I. Satake, On numerical invariants of arithmetic varieties
  of ${\Bbb Q}$-rank one, in: Automorphic forms of several variables,
  Taniguchi Symposium, Katata, 1983. Birkh\"auser: Boston 1984.

\bibitem[S2]{S2} I. Satake, On some properties of holomorphic imbeddings of
  symmetric domains, Amer. J. Math. {\bf 91} (1969), 289-305.

\bibitem[Sch]{Sch} W. Scharlau, ``Quadratic and hermitian forms'',
  Grundlehren der Math. Wiss. {\bf 270}, Springer: Berlin 1988.

\bibitem[Sh]{Sh2} G. Shimura, On analytic families of polarised abelian
  varieties and automorphic functions, Ann.~of Math. {\bf 78}
  (1963),149-192.

\bibitem[SC]{SC} A. Ash, D. Mumford, M. Rappoport, Y. Tai, ``Smooth
  compactification of locally symmetric spaces,'' Math.~Sci.~Press:
  Brookline 1975.

\bibitem[TV]{T} E. Thomas \& A. Vasquez, Chern numbers of Hilbert modular
  varieties, J.~r.~a.~Math, {\bf 324} (1981), 192-210.

\bibitem[T]{tits} J. Tits, Classification of simple algebraic groups, in:
  Algebraic groups and discontinuous subgroups, Proc.~Symp.~Pure Math.
  Vol.~{\bf 9}, AMS: Providence 1966.

\end{thebibliography}
\end{document}